\documentclass[12pt]{article}

\usepackage{graphicx}
\usepackage{amsmath,amssymb}
\usepackage{epstopdf}
\usepackage{psfrag}
\usepackage{float}
\usepackage[small,nohug,heads=littlevee]{diagrams} \diagramstyle[labelstyle=\scriptstyle] 
\newtheorem{theorem}{Theorem}
\newtheorem{corollary}[theorem]{Corollary}

\numberwithin{equation}{section}

\makeatletter 
\long\def\@makecaption#1#2{%
\vskip\abovecaptionskip 
\sbox\@tempboxa{\textrm{#1}. \textrm{#2}} 
\ifdim \wd\@tempboxa >\hsize 
{\textrm{#1}. \textrm{#2}\par} 
\else 
\global \@minipagefalse 
\hb@xt@\hsize{\hfil\box\@tempboxa\hfil}%
\fi 
\vskip\belowcaptionskip} 
\makeatother 

\newcommand{\sn}{\mbox{sn}}
\newcommand{\cn}{\mbox{cn}}
\newcommand{\dn}{\mbox{dn}}

\begin{document}
\date{}

\begin{flushright}{YITP-SB-07-41}\end{flushright}

\begin{center}
{\LARGE Elliptic constructions of hyperk\"ahler metrics III: \\[5pt] Gravitons and Poncelet polygons}
\vskip40pt
{\bf Radu A. Iona\c{s}}
\end{center}
\vskip10pt
\centerline{\it C.N.Yang Institute for Theoretical Physics, Stony Brook University}
\centerline{\it Stony Brook, NY 11794-3840, USA }
\centerline{\tt ionas@max2.physics.sunysb.edu}
\vskip60pt

\begin{abstract}
\noindent In the generalized Legendre approach, the equation describing an asymptotically locally Euclidean space of type $D_n$ is found to admit an algebraic formulation in terms of the group law on a Weierstrass cubic. This curve has the structure of a Cayley cubic for a pencil generated by two  transversal plane conics, that is, it takes the form $Y^2 = \det ({\cal A}+X{\cal B})$, where ${\cal A}$ and ${\cal B}$ are the defining $3 \times 3$ matrices of the conics. In this light, the equation  can be interpreted as the closure condition for an elliptic billiard trajectory tangent to the conic ${\cal B}$ and bouncing into various conics of the pencil determined by the positions of the monopoles. Poncelet's porism guarantees then that once a trajectory closes to a star polygon, any trajectory will close, regardless of the starting point and after the same number of steps.
\end{abstract}
\vspace*{40pt}

\newpage

\tableofcontents

\setcounter{section}{-1}

\section{Introduction}

Asymptotically locally Euclidean (ALE) spaces are non-compact complete Riemannian hyperk\"ahler 4-manifolds whose boundary at infinity resembles a quotient $\mathbb{H}/\Gamma$ of the Euclidean space $\mathbb{H} \simeq \mathbb{R}^4$ by a finite subgroup $\Gamma \subset SU(2) \simeq Sp(1,\mathbb{H})$. The Riemannian metric is required to approximate the Euclidean metric up to order $1/r^4$.  The spaces corresponding to $\Gamma = \mathbb{Z}_k$ were found by Gibbons and Hawking in 1978 \cite{Gibbons:1979zt}. In 1987, Kronheimer \cite{MR992334} gave an explicit construction for all $\Gamma$ which also clarified the conjectured ADE classification of these spaces \cite{MR520463}. The construction relied on the hyperk\"ahler quotient technique, an extension of the symplectic quotient of Marsden and Weinstein to holomorphic settings \cite{Lindstrom:1983rt,Hitchin:1986ea}. One takes the quotient of a parent space, an $\mathbb{H}$-linear flat hyperk\"ahler manifold set up in such a way as to possess a natural action of $\Gamma \subset SU(2)$, with the tri-holomorphic action of a maximal product of unitary groups commuting with the $\Gamma$-action. The various elements of this approach can be conveniently summarized in a quiver diagram, related to the affine Dynkin diagram of the simply-laced group corresponding to $\Gamma$ in the McKay classification. The intersection matrix of the non-trivial 2-cycles of the quotient manifold is given by the negative of the Cartan matrix of this group's Lie algebra.

The ALE as well as the associated ALF spaces of type $D_n$ have been approached within the framework of the generalized Legendre transform of Lindstr\"om and Ro\v{c}ek \cite{Lindstrom:1983rt,Hitchin:1986ea,Lindstrom:1987ks} in \cite{Ivanov:1995cy,Chalmers:1998pu}. The form of the corresponding generalized Legendre transform holomorphic potentials which encode all the metric information was conjectured based on the known \cite{Hitchin:1986ea} holomorphic potentials of the ALE respectively ALF spaces of type $A_n$ that the type $D_n$ spaces approach asymptotically. This conjecture was confirmed by Cherkis and Kapustin in \cite{Cherkis:1998xca}. The ALE spaces were identified by these authors with moduli spaces of Nahm equations \cite{MR807414} by means of string-theoretical arguments. In particular, an algebraic constraint of Ercolani-Sinha type \cite{Ercolani:1989tp} on the Jacobian of a spectral curve emerged from their analysis. A further confirmation came in \cite{Cherkis:2003wk}, this time from the alternate direction of Hitchin's twistor approach to monopoles \cite{MR649818,MR709461}.

In this paper we obtain a refined form of the Ercolani-Sinha constraint from a generalized Legendre equation. We also find that the 2-monopole spectral curve can be biholomorphically mapped to a Cayley cubic form and thus related to a pencil generated by two transversal plane conics. These two results prompt us to interpret the Ercolani-Sinha constraint as the closure condition for an in-and-circumscribed Poncelet star polygon.

Not long ago, Poncelet polygons have been used by Hitchin to derive a special class of solutions to a certain Painlev\'e VI equation \cite{MR1351506,MR1989487}. Although otherwise very different, this problem has one thing in common with the one that we address, namely the presence of a dihedral symmetry. Is this just a coincidence, or part of a larger pattern? This interesting question remains yet to be answered.

The paper is divided into three distinct parts. In the first part (section \ref{Leg-add}) we discuss an addition theorem of Legendre concerning incomplete elliptic integrals and spherical triangles. In the second part (section \ref{SEC:Poncelet-porism}) we review, following a series of papers by Cayley \cite{Cayley:1853a,Cayley:1853b, Cayley:1853c, Cayley:1861} and their modern algebraic-geometric translation due to Griffiths and Harris \cite{MR0498606,MR497281}, the Poncelet problem of in-and-circumscribing a polygon to a pair of conics. In the third part (sections \ref{SEC:O4-multiplets} and \ref{SEC:ALE-typeD}), we continue the analysis of the 2-monopole spectral curve started in \cite{MM1} and study the ALE spaces of type $D_n$ using the generalized Legendre transform approach. The relevance of our incursions into the realm of the classical projective geometry problems of Legendre and Poncelet will eventually come to light in this context.

\section{Incomplete elliptic integrals of first kind and spherical triangles} \label{Leg-add}

In this section we will be concerned with the following theorem, due to Legendre \cite{Legendre:1825}:

\begin{theorem}{\rm (Legendre)} 
The equality
\begin{equation}
F(\sin A,k)+F(\sin B,k)+F(\sin C,k) = 2 K(k) \label{Legendre_add_thm}
\end{equation}
holds if and only if the amplitudes $A$, $B$ and $C$ form the angles of a spherical triangle, the lengths of the sides of which can be determined from the sine theorem
\begin{equation}
\frac{\sin a}{\sin A} = \frac{\sin b}{\sin B} = \frac{\sin c}{\sin C} = k \label{sine_thm}
\end{equation}
where $k \in [0,1]$ is the elliptic modulus.\footnote{This addition theorem appears in Legendre's treatise in a slightly different form and states that if the sum of two elliptic integrals equals a third one then their amplitudes must form the {\it sides} of a spherical triangle. Their asymmetric occurence in the elliptic formula forces then one to impose certain sign choice prescriptions. The form presented here, manifestly symmetric, avoids this formal inconvenience.}
\end{theorem}


\begin{figure}[bth]
\centering
\scalebox{0.4}{\includegraphics{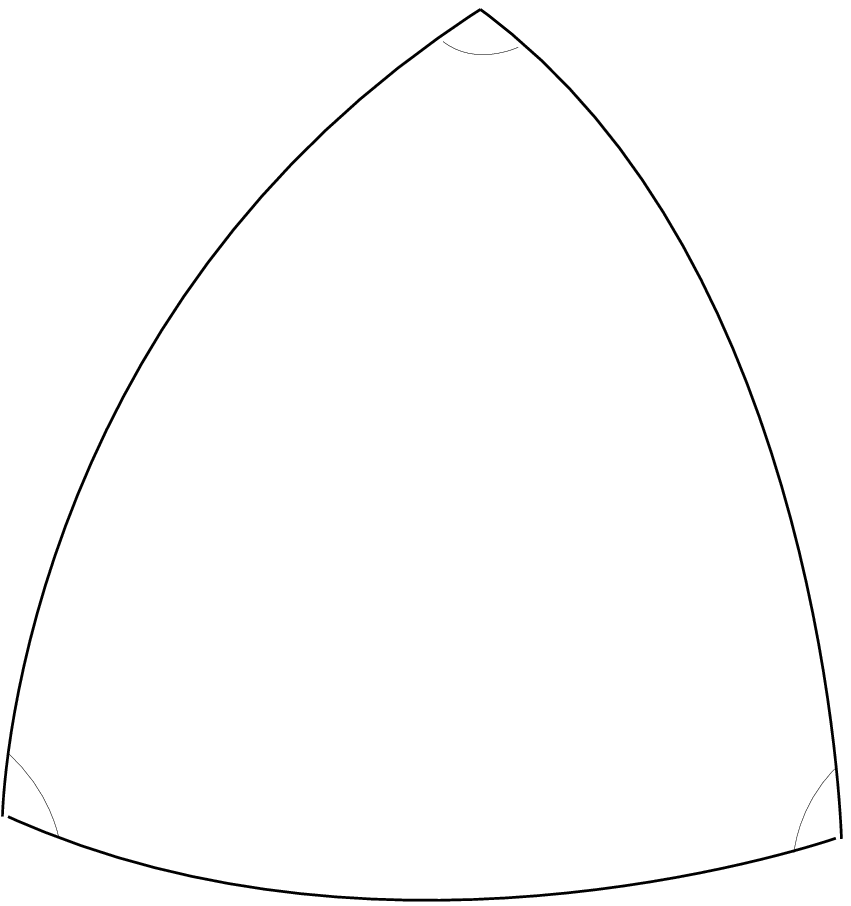}}
\put(1,-1){$C$}
\put(-46,108){$A$}
\put(-111,3){$B$}
\put(-50,-9){$a$}
\put(-6,58){$b$}
\put(-91,63){$c$}
\end{figure} 

Legendre's addition theorem is equivalent to the addition theorems for Jacobi's elliptic functions in the real domain. The theory of Jacobi elliptic functions emerged historically from the study of the problem of inverting incomplete Legendre elliptic integrals of the first kind, {\it i.e.}, 
\begin{equation}
\sn(F(u,k),k) = u
\end{equation}
An excellent reference for this topic is Cayley's treatise on elliptic functions \cite{MR0124532}.

The list of properties of Jacobi's elliptic functions $\sn(u,k)$, $\cn(u,k)$ and $\dn(u,k)$ includes \\
$\bullet$ trigonometric-like relations
\begin{equation}
\sn^2 u + \cn^2 u =1 \hspace{40pt} \dn^2 u + k^2 \sn^2 u =1
\end{equation}
$\bullet$ reflection symmetry 
\vspace{-10pt}
\begin{eqnarray}
\sn (-u) &\hspace{-7pt} = \hspace{-7pt}& -\sn\, u \nonumber \\
\cn (-u) &\hspace{-7pt} = \hspace{-7pt}& +\cn\, u  \label{refl_symm} \\
\dn (-u) &\hspace{-7pt} = \hspace{-7pt}& +\dn\, u \nonumber
\end{eqnarray}
$\bullet$ double-periodicity 
\begin{equation}
\begin{array}{lll}
\sn (\upsilon +2mK+2m'iK') &\hspace{-10pt} =  (-)^m 	&\hspace{-10pt} \sn\, \upsilon \\ [4pt]
\cn (\upsilon +2mK+2m'iK') &\hspace{-10pt} =  (-)^{m+m'} 	&\hspace{-10pt} \cn\, \upsilon \\ [4pt]
\dn (\upsilon +2mK+2m'iK') &\hspace{-10pt} =  (-)^{m'} 	&\hspace{-10pt} \dn\, \upsilon
\end{array} \label{periodicity}
\end{equation}
$\bullet$ addition theorems
\begin{eqnarray}
\sn (u+v) &\hspace{-7pt} = \hspace{-7pt}& \frac{\sn\, u\, \cn\, v\, \dn\, v + \sn\, v\, \cn\, u\, \dn\, u}{1-k^2\sn^2 u\, \sn^2 v} \label{sn(u+v)} \\ [2pt]
\cn (u+v) &\hspace{-7pt} = \hspace{-7pt}& \frac{\cn\, u\, \cn\, v - \sn\, u\, \dn\, u\, \sn\, v\, \dn\, v}{1-k^2\sn^2 u\, \sn^2 v}  \label{sncndn_add_thms} \\ [2pt]
\dn (u+v) &\hspace{-7pt} = \hspace{-7pt}& \frac{\dn\, u\, \dn\, v - k^2 \sn\, u\, \cn\, u\, \sn\, v\, \cn\, v}{1-k^2\sn^2 u\, \sn^2 v} 
\end{eqnarray}
To prove Legendre's addition theorem, one can take for instance to the other side of the equality in (\ref{Legendre_add_thm}) one of the incomplete elliptic integrals, say $F(\sin C,k)$, and then act in turns on the equation thus obtained with Jacobi's elliptic functions $\cn$ and $\dn$ of modulus $k$. This allows one to employ on the l.h.s. the  addition formulas (\ref{sncndn_add_thms})  and on the r.h.s. the properties (\ref{refl_symm}) and (\ref{periodicity}) of Jacobi's elliptic functions, together with the relations\begin{eqnarray*}
\sn\, F(\sin A,k) &\hspace{-7pt} = \hspace{-7pt}& \sin A \\ [2pt]
\cn\, F(\sin A,k) &\hspace{-7pt} = \hspace{-7pt}& \cos A \\ [2pt]
\dn\, F(\sin A,k) &\hspace{-7pt} = \hspace{-7pt}& \cos a
\end{eqnarray*}
and the similar ones for $F(\sin B,k)$. The remaining $k^2$'s should be substituted by $\sin a \sin b \div  \sin A \sin B$, in accordance with equations~(\ref{sine_thm}). Eventually, one arrives, after some algebraic manipulations, at the cosine theorems of spherical trigonometry. Conversely, it is also possible to derive Jacobi's addition theorems starting from Legendre's, see {\it e.g.} \cite{MR0124532}. The two are thus equivalent in the real domain.

\section{Poncelet's closure theorem} \label{SEC:Poncelet-porism}

\subsection{Poncelet's porism}

In this section we will be concerned with the following theorem, due to Poncelet:

\begin{theorem}{\hspace{-4pt}{\rm (Poncelet's porism\footnote{"A proposition affirming the possibility of finding such conditions as will render a certain problem indeterminate, or capable of innumerable solutions" \cite{Playfair:1792}})}}
Given two plane conics ${\cal A}$ and ${\cal B}$, with ${\cal A}$ lying inside ${\cal B}$, if there exists a (possibly star) polygon inscribed in ${\cal B}$ and circumscribed about ${\cal A}$ then there exist an infinity of such polygons. 
\end{theorem}

\begin{figure}[tbh]
\centering
\scalebox{0.6}{\includegraphics{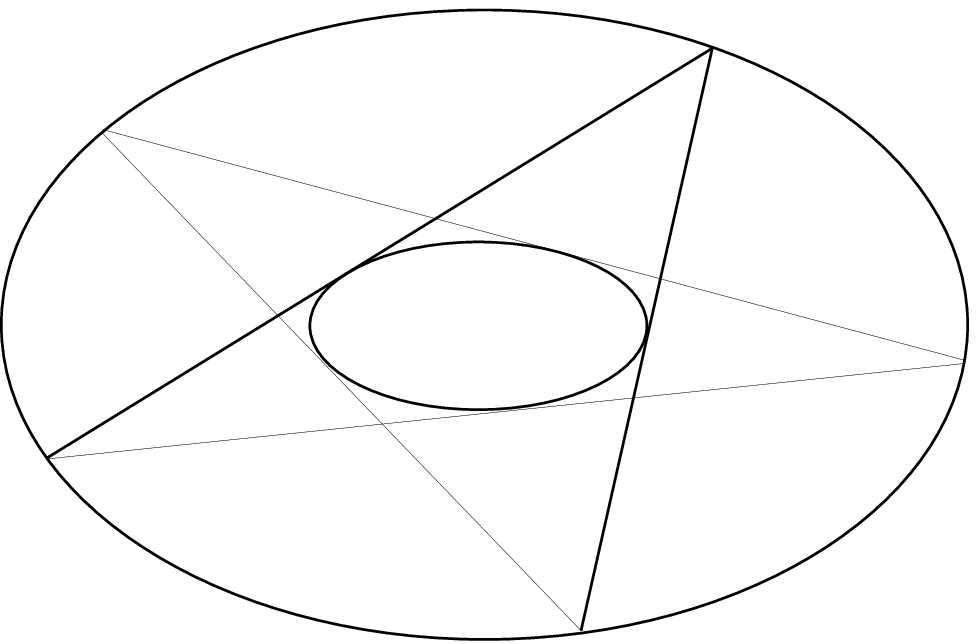}} 
\put(-10,87){${\cal B}$}
\put(-75,55){${\cal A}$}
\caption{}
\label{Poncelet0}
\end{figure}

Equivalently, consider a point $P_0$ lying on the conic ${\cal B}$ and from it draw a tangent $L_1$ to the conic ${\cal C}$, which will intersect again ${\cal B}$ at a point $P_1$. Repeat this construction starting this time from $P_1$, a.s.o. This yields a series of pairs of tangents to ${\cal A}$ and points on ${\cal B}$, $(P_1,L_1)$, $(P_2,L_2)$, $\cdots$. If, after a finite number of steps, one arrives back at $P_0$, then Poncelet's porism states that this will happen regardless of which starting point $P_0$ one chooses.

For original papers, reviews and related material, see \cite{Cayley:1853a,Cayley:1853b, Cayley:1853c, Cayley:1861, Jacobi:1881, MR0498606, MR497281, MR917349,MR946326, MR1218986, MR1218987, MR1989487, MR1351506,Tabachnikov2005}.

\subsection{The projective geometry of plane conics}

A point in the projective plane $\mathbb{P}^2$  is specified by its homogeneous coordinates, $x = [x_1:x_2:x_3]$. The equation of a projective line in $\mathbb{P}^2$ that passes through the point $x$ is
\begin{equation}
(y,x) = \sum_{i=1}^3 y_i x_i = 0 \label{pline}
\end{equation}
A line is specified by its coefficients, the tangential coordinates $y^* = (y_1:y_2:y_3)$, which can be thought of  as being the homogeneous coordinates of a point in the dual projective plane $\mathbb{P}^{2*}$. One can similarly argue that a line in $\mathbb{P}^{2*}$ corresponds to a point in $\mathbb{P}^2$. The symmetry of the equation (\ref{pline}) at the interchange of $x$ and $y$ results in an ambiguity of interpretation of what one means by "points" and "lines" which lays at the heart of the principle of duality of projective geometry. For example, a fundamental theorem of projective geometry states that through any two distinct points in a projective plane there passes exactly one line. Applying it to the dual projective plane yields the dual theorem: any two distinct lines in a projective plane intersect exactly once. The duality correspondence preserves incidence relationships.

Projective conics in $\mathbb{P}^2$ are described by means of quadratic equations
\begin{equation}
(x,Qx) = \sum_{i,j=1}^3 x_i Q_{ij} x_j = 0
\end{equation}
where $Q$ is a symmetric $3 \times 3$ matrix. Such a projective variety is a smooth submanifold of $\mathbb{P}^2$ and thus a Riemann surface if and only if the matrix $Q$ is non-singular. 

\begin{figure}[tbh]
\centering
\scalebox{0.6}{\includegraphics{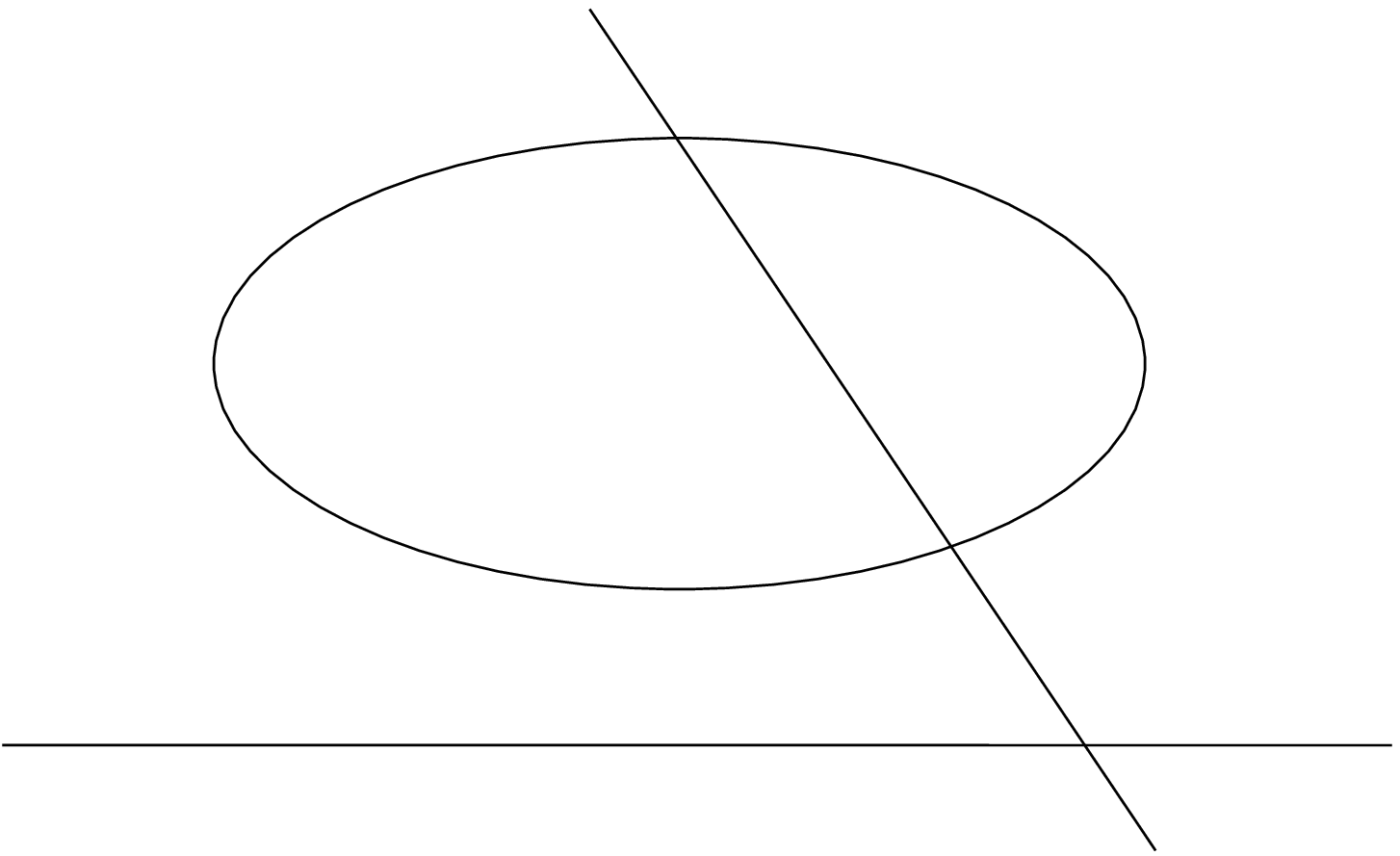}} 
\put(-130,135){$P_0$}
\put(-200,25){$L_0$}
\caption{}
\label{Poncelet3}
\end{figure}

Given any smooth conic ${\cal C}$, let $P_0$ be a point on ${\cal C}$  and $L_0$ be a line that does not contain $P_0$, see \figurename~\ref{Poncelet3}. By B\'{e}zout's theorem, any line in the projective plane intersects a smooth conic exactly twice, counting multiplicities. Then any line that passes through $P_0$ will intersect the conic at one other point which is in one-to-one correspondence with the point at which the line intersects $L_0$. This stereographic projection-like construction establishes a biholomorphic mapping ${\cal C} \longrightarrow L_0 \simeq \mathbb{P}^1$, {\it i.e.}, it provides a {\it rational parametrization} of the conic. 

The equation of the tangent line to a smooth conic ${\cal C}$ at a point $x \in {\cal C}$ is
\begin{equation}
0 = \frac{1}{2} \sum_{i=1}^3 y_i \frac{\partial}{\partial x^i} (x, Qx) = (Qx,y)
\end{equation} 
The tangent is thus the subspace of $\mathbb{P}^2$ orthogonal to $x$ with respect to the symmetric bilinear form associated to ${\cal C}$. The dual coordinates of the tangent line are $x^* = Qx$. They satisfy
\begin{equation}
(x^*, Q^{-1} x^*) = 0
\end{equation} 
so they are points on the dual conic ${\cal C}^*$ defined by the inverse matrix $Q^{-1}$. So, the envelope of tangents to a conic is also a conic. Points are dual to lines, conics are self-dual.

\subsection{Poncelet's construction, {\it \`a la} Griffiths and Harris}

Consider two smooth conics defined by the matrices\footnote{In the following, we will refer to a conic using the symbol of its defining matrix.} ${\cal A}$ and ${\cal B}$, with no common components, so that they intersect transversally at four points in general position. The set of plane conics that contain these four points, {\it i.e.}, the {\it pencil} of conics generated by ${\cal A}$ and ${\cal B}$, is given by the one-parameter family 
\begin{equation}
{\cal C}_{X} = {\cal A} + X  {\cal B} 
\end{equation}
with $X \in \mathbb{C} \cup \{ \infty \} \simeq \mathbb{P}^1$. In particular, ${\cal C}_0 = {\cal A}$ and ${\cal C}_{\infty} = {\cal B}$. Among the conics in the pencil there are three singular ones, consisting of the three pairs of lines obtained by joining in all possible ways pairs of the four intersection points.

Fix a conic ${\cal C}_{X_0}$ from the pencil, non-singular and different from ${\cal B}$.  In order to address the Poncelet problem, Griffiths and Harris \cite{MR0498606} construct the incidence correspondence\footnote{The construction that we present here is in fact dual to that of Griffiths and Harris.}
\begin{equation}
\Sigma = \{(P,L) \in {\cal B} \times {\cal C}_{X_0}^* \ | \ P \in L\}
\end{equation}
{\it i.e.}, the set of pairs of points $P$ on ${\cal B}$ and tangents $L$ to ${\cal C}_{X_0}$ subject to the incidence condition that $L$ passes through $P$. As both conics ${\cal B}$ and ${\cal C}_{X_0}^*$ can be rationally parametrized, $\Sigma \subset {\cal B} \times {\cal C}_{X_0}^* \simeq \mathbb{P}^1 \times \mathbb{P}^1$. The transversality of the intersection ${\cal B} \cap {\cal C}_{X_0} = {\cal B} \cap {\cal A}$ insures that $\Sigma$ is a non-singular variety. 

\begin{figure}[tbh]
\centering
\scalebox{0.6}{\includegraphics{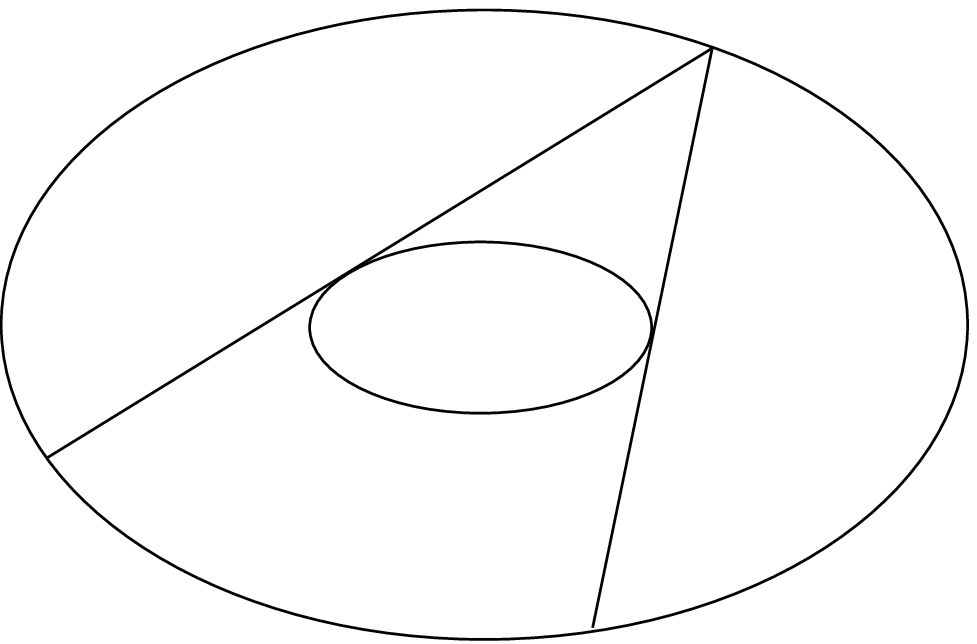}} 
\put(-46,106){$P$}
\put(-74,-10){$P'$}
\put(-120,67){$L$}
\put(-52,48){$L'$}
\put(-100,28){${\cal C}_{X_0}$}
\put(1,50){${\cal B}$}
\caption{}
\label{Poncelet1a}
\end{figure}

Given a point and a smooth conic, there exist exactly two lines, counting multiplicities, that are tangent to the conic and intersect each other at the given point, with the two tangents being confounded if and only if the point belongs to the conic. This is dual to the statement that a line intersects a smooth conic exactly twice, counting multiplicities, with the intersection points coinciding if and only if the line is tangent to the conic.
As a consequence, the variety $\Sigma$ has two natural involutive automorphisms, namely
\begin{eqnarray}
 i_1(P,L') &\hspace{-7pt}=\hspace{-7pt}& (P',L') \\ [2pt]
 i_2(P,L\phantom{'}) &\hspace{-7pt}=\hspace{-7pt}& (P\phantom{'},L') 
\end{eqnarray}
for notations see \figurename~\ref{Poncelet1a}.
The fixed points of $i_1$ are the four points of ${\cal B}^* \cap {\cal C}_{X_0}^*$, {\it i.e.}, the four common tangents of ${\cal B}$ and ${\cal C}_{X_0}$, whereas the fixed points of $i_2$ are the four points of ${\cal B} \cap {\cal C}_{X_0} = {\cal B} \cap {\cal A}$. The relevance to the Poncelet problem becomes transparent when we observe that the action of the composed automorphism $j = i_1 \circ i_2$, namely,
\begin{eqnarray}
j(P,L) = (P',L')
\end{eqnarray}
offers a realization of the basic step in the geometric construction of Poncelet polygons. 

The projection
\begin{eqnarray}
\begin{diagram}[width=30pt]
\Sigma & (P,L) \\
\dTo_{2:1} & \dMapsto \\
{\cal B} & P
\end{diagram} \label{double_cover}
\end{eqnarray}
exhibits $\Sigma$ as a branched double-cover of ${\cal B} \simeq \mathbb{P}^1$. The action of $i_2$ interchanges the sheets of the double-cover, the branching points being the fixed points of $i_2$. The Riemann-Hurwitz formula tells us then that $\Sigma$ has genus $1$, {\it i.e.}, $\Sigma$ is an elliptic curve.

\begin{figure}[tbh]
\centering
\scalebox{0.5}{\includegraphics{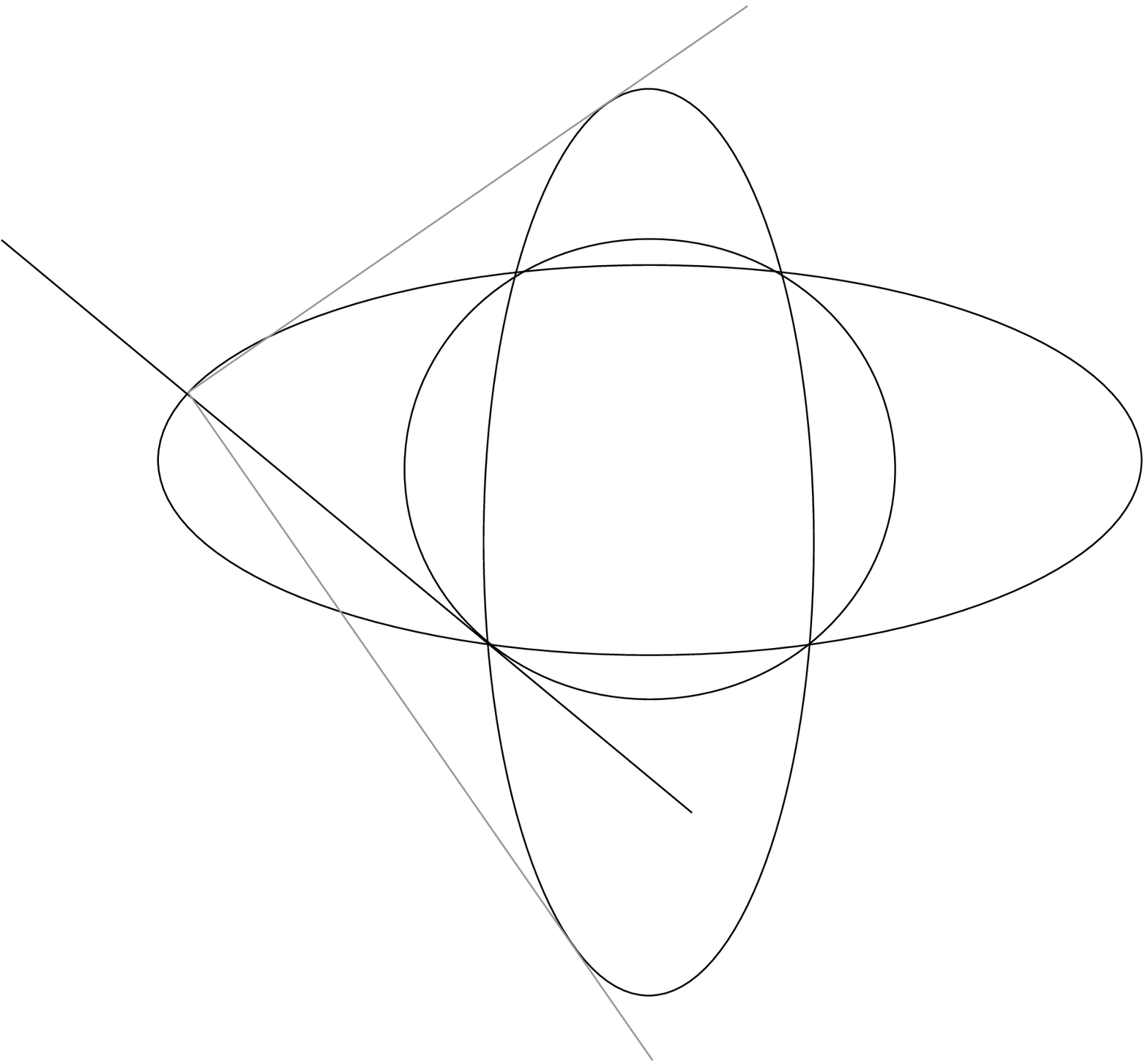}} 
\put(-167,88){$P_{e_0}$}
\put(-160,188){$P_{e_1}$}
\put(-84,187){$P_{e_3}$}
\put(-77,87){$P_{e_4}$}
\put(-229,163){$P_{X}$}
\put(-187,137){${\cal C}_{X}$}
\put(2,137){${\cal C}_{\infty} = {\cal B}$}
\put(-148,137){${\cal C}_{X_0}$}
\caption{}
\label{Poncelet2}
\end{figure}

To cast $\Sigma$ in a more explicit form, Griffiths and Harris, following the ideas of Cayley \cite{Cayley:1853a}, use an ingenious rational parametrization construction for ${\cal B}$.
Choose one of the four intersection points, and take the tangent to an arbitrary conic ${\cal C}_{X}$ from the pencil through this point, see \figurename~\ref{Poncelet2}. The tangent intersects ${\cal B}$ at one more point, which we then label $P_{X}$. Together with the tangent to ${\cal B}$ through one of the other intersection points, this gives a rational parametrization of ${\cal B}$ by the complex parameter $X$. 

The chosen point of intersection itself is the limit case given by the tangent to ${\cal C}_{\infty} = {\cal B}$, so it corresponds to $P_{e_0 = \infty}$. Let $P_{e_1}$, $P_{e_2}$, $P_{e_3}$ be the other three intersection points. The conics ${\cal C}_{e_1}$, ${\cal C}_{e_2}$, ${\cal C}_{e_3}$ that parametrize them are the three singular conics of the pencil. Indeed, a conic corresponding to one of the points $P_{e_i}$ with $i = 1, 2$ or $3$ has to simultaneously satisfy the following two properties: $1)$ the line joining $P_{e_0}$ and $P_{e_i}$ is tangent to it and $2)$ the two, by assumption disjoint, intersection points $P_{e_0}$ and $P_{e_i}$ belong to it, as they belong to all conics in the pencil. But these two requirements cannot be satisfied at once unless the conic is degenerate. This occurs when the defining matrix is singular, {\it i.e.}, when $X = e_i$ is a solution of $\det ({\cal A}+X {\cal B}) = 0$. 

A point $P_{X} \in {\cal B}$ together with a choice of tangent to ${\cal C}_{X_0}$ define a point on the double cover $\Sigma$. Branching occurs when the tangents through $P_{X}$ coincide, and this cannot happen unless $P_{X} \in {\cal C}_{X_0}$, in which case $P_{X} \in {\cal C}_{X_0} \cap {\cal B} = {\cal A} \cap {\cal B}$, and so, by the argument above, $X$ has to be a solution of $\det ({\cal A}+X {\cal B}) = 0$. From these considerations if follows that the elliptic curve $\Sigma$ is isomorphic to Cayley's cubic\footnote{For this reason, in what follows we will denote Griffiths and Harris's incidence correspondence and Cayley's cubic curve by the same symbol, $\Sigma$.}
\begin{equation}
Y^2 = \det ({\cal A}+X {\cal B})
\end{equation}
the isomorphism between them being given by 
\begin{eqnarray}
(P_{X},L\phantom{'}) &\longleftrightarrow& (X,+Y) \nonumber \\ [2pt]
(P_{X},L') &\longleftrightarrow& (X,-Y)
\end{eqnarray}

As an elliptic curve, $\Sigma$ posesses an abelian differential, {\it i.e.}, a globally holomorphic 1-form 
\begin{equation}
\varpi = \frac{dX}{Y}
\end{equation}
with associated period lattice $X = \mathbb{Z}\cdot 2\omega + \mathbb{Z}\cdot 2\omega'$. The fundamental periods $2\omega$ and $2\omega'$, chosen such that $\mbox{Im}\, \omega'/\omega > 0$, are the integrals of $\varpi$ over the $a$ and $b$-cycles of the torus $\Sigma$. One can exploit $\varpi$ to give an alternative description of $\Sigma$ by means of the Abel-Jacobi map, an analytic isomorphism between $\Sigma$ and its Jacobian variety, $\mathbb{C}/\Lambda$,
\begin{eqnarray}
\begin{diagram}[width=35pt,height=28pt]
\Sigma & \rTo^{\simeq} & \mathbb{C}/\Lambda \\ 
(X,Y) & \rMapsto & \displaystyle{\int_{\Gamma} \varpi}
\end{diagram}
\end{eqnarray}
The integral, taken on a path $\Gamma$ on $\Sigma$ based at an arbitrary fixed point, is independent of the path modulo integer multiples of the periods, that is to say, it defines an equivalence class on $\mathbb{C}/\Lambda$.

The automorphisms of $\Sigma$ are carried over by the Abel-Jacobi map to $\mathbb{C}/\Lambda$, and for simplicity we will denote the corresponding automorphisms of $\mathbb{C}/\Lambda$ by the same letters. Any automorphism $i$ of $\mathbb{C}/\Lambda$ is induced by an automorphism of its universal cover $\mathbb{C}$,  $ \tilde{\imath}(u) = a u +b$, for any $u \in \mathbb{C}$. Then $i$ is involutive, {\it i.e.}, $i^2(u) = u\bmod{\Lambda}$ if and only if $a^2 = 1$ and $(a+1)b = 0 \bmod{\Lambda}$. In the case when $a=+1$, one easily argues that $i$ has no fixed points unless it is the trivial automorphism of $\Sigma$, in which case all points of $\Sigma$ are fixed points. Since we want $i_1$ and $i_2$ to have no more and no less than four fixed points each, this cannot be the case. So $a = -1$ for both, that is, $i_1(u) = -u + b_1 \bmod{\Lambda}$ and $i_2(u) = -u + b_2 \bmod{\Lambda}$. Moreover, one can always redefine $u$ by a shift to put $b_2 = 0$. Eventually, renaming $b_1=u_0$, we have
\begin{eqnarray}
i_1(u) &\hspace{-7pt} = \hspace{-7pt}& -u+u_0 \bmod{\Lambda} \label{i1_Jac}\\ [2pt]
i_2(u) &\hspace{-7pt} = \hspace{-7pt}& -u \bmod{\Lambda} \label{i2_Jac}
\end{eqnarray}
and so 
\begin{equation}
j(u) = u+u_0 \bmod{\Lambda} \label{j_Jac}
\end{equation}
Note that on $\mathbb{C}/\Lambda$ $i_1$ has the four fixed points $u_0/2 \bmod{\Lambda/2}$ and $i_2$ the four fixed points $0 \bmod{\Lambda/2}$.

On the other hand, the inverse of the Abel-Jacobi map 
\begin{eqnarray}
\begin{diagram}[width=30pt]
\mathbb{P}^2 \supset \Sigma               	& \lTo^{\simeq}  	& \mathbb{C}/\Lambda \\
(X,Y)=(\wp(u),\wp'(u)) 	& \lMapsto          	& u					
\end{diagram}
\end{eqnarray}
gives  $X = \wp(u)$, $Y = \wp'(u)$ and thus $\varpi = du$, with $\wp(u)$ and its derivative $\wp'(u)$ elliptic functions of order $2$ respectively $3$, doubly-periodic with period lattice $\Lambda$, meromorphic on $\mathbb{C}$. 
Based on (\ref{i2_Jac}) and on the fact that the induced action of $i_2$ on the cubic curve interchanges $(X,+Y)$ with $(X,-Y)$, one can argue that $\wp(u)$ is even and $\wp'(u)$ is odd. The map
\begin{eqnarray}
\begin{diagram}[width=36pt]
\mathbb{C}/\Lambda & u \bmod{\Lambda} \\
\dTo_{2:1} & \dMapsto \\
\mathbb{P}^1& \wp(u)
\end{diagram}
\end{eqnarray}
is a branched double-covering of $\mathbb{P}^1$ by $\mathbb{C}/\Lambda$, with $\pm u \bmod{\Lambda}$ mapped to $\wp(u) = \wp(-u)$ and the four fixed points of $i_2$ on $\mathbb{C}/\Lambda$, {\it i.e.}, $0 \bmod{\Lambda/2}$, mapped to the branching points $e_0$, $e_1$, $e_2$, $e_3$.
We will assume that $\wp(0) = e_0$.

The various correspondences are summarized in the diagram
\begin{eqnarray}
\begin{diagram}[midshaft,width=32pt]
(X,Y) = (\wp(u),\wp'(u))	\hspace{-75pt}		& \hfill\hfill\hfill\hspace{-10pt}\lMapsto 	&&& u \bmod{\Lambda} \\
\dMapsto & \Sigma              	& \rTo^{\simeq}_{\scriptstyle{\rm Abel-Jacobi}}  	& \mathbb{C}/\Lambda 	& \dMapsto \\
								& \dTo^{2:1} 	&		& \dTo_{2:1} \\
								& {\cal B}	& \rTo^{\simeq}_{\scriptstyle{\rm rational\ parametrization}}	& \mathbb{P}^1			&  \\
P_{X}			&			\rMapsto	& 					&& X= \wp(u)
\end{diagram}
\end{eqnarray}

\begin{figure}[tbh]
\centering
\scalebox{0.6}{\includegraphics{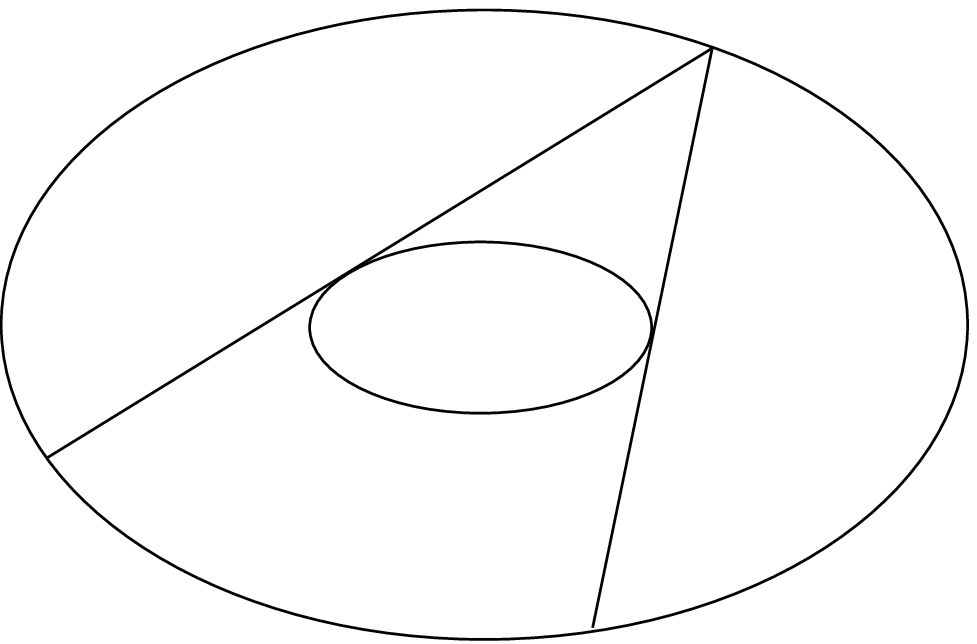}} 
\put(-47,109){$P_{\wp(u)}$}
\put(-70,-9){$P_{\wp(u+u_0)}$}
\put(-190,20){$P_{\wp(u-u_0)}$}
\put(-107,30){${\cal C}_{\wp(u_0)}$}
\put(2,60){${\cal B}$}
\caption{}
\label{Poncelet1b}
\end{figure}

A line that passes through the point $P_{e_0} \in {\cal B}$ and is tangent to the fixed conic ${\cal C}_{X_0}$ will, by the above choice of rational parametrization for ${\cal B}$, intersect ${\cal B}$ again at the point parametrized by $X_0$, {\it i.e.}, $P_{X_0}$. From $e_0 = \wp(0)$ together with (\ref{i1_Jac}) we then obtain the interpretation of the $u_0$-shift of $j$, namely
\begin{equation}
\wp(u_0) = X_0
\end{equation}
One has then the following (see also \figurename~\ref{Poncelet1b})
\begin{corollary}{\hspace{-9pt} \rm (Cayley)} 
The tangents from a point $P_{\wp(u)} \in {\cal B}$ to a non-singular conic ${\cal C}_{\wp(u_0)}$ from the pencil generated by the conics ${\cal A}$ and ${\cal B}$  will intersect again ${\cal B}$ at the points $P_{\wp(u\pm u_0)}$.
\end{corollary}

As observed above, the basic step in the construction of Poncelet polygons corresponds on the Jacobian variety of $\Sigma$ to the action of the automorphism $j$. The Poncelet problem can be reformulated in the following terms: the polygon closes after $n$ steps if $j^n$ has fixed points on $\Sigma$.  From $j^n(u) = u+nu_0$ it follows that the necessary and sufficient condition that $j^n$ has a fixed point is 
\begin{equation}
n u_0 = 0 \bmod{\Lambda} \label{cycl-n}
\end{equation}
The elements of the Jacobian form an abelian group with respect to addition modulo lattice shifts. The condition (\ref{cycl-n}) means that $u_0$  is a cyclic element of this group of order $n$.  This condition is clearly independent of the point $u \in \mathbb{C}/\Lambda \longleftrightarrow (P,L) \in \Sigma$, and this proves the porism. 

The Poncelet problem can be generalized in the following way: consider a conic ${\cal B}$ and a series of  conics ${\cal C}_1$, ${\cal C}_2$, ${\cal C}_3$, $\cdots$ from the pencil generated by ${\cal B}$ and another transversal conic, ${\cal A}$. Take a point $P_0$ on ${\cal B}$ and draw a tangent  $L_1$ to the conic ${\cal C}_1$ which intersects again ${\cal B}$ at the point $P_1$. From $P_1$ draw a tangent $L_2$ to the conic ${\cal C}_2$, a.s.o. Dually, this reads as follows: take a point $L_1$ on ${\cal C}_1^*$ and draw a tangent $P_1$ to ${\cal B}^*$ which intersects ${\cal C}_2^*$ at a point $L_2$. From $L_2$ draw a tangent to ${\cal B}^*$ that intersects ${\cal C}_3^*$ at the point $L_3$, a.s.o. In this case one has not one but a series of automorphisms of type $j$, one for each conic ${\cal C}_i$. The above arguments can be easily extended to give the condition for this construction to close after $n$ steps: one has to have
\begin{equation}
u_1 + \cdots + u_n = 0 \bmod{\Lambda} \label{gen-cycl-n}
\end{equation}
where $u_1$, $\cdots$, $u_n$ are such that ${\cal C}_i = {\cal C}_{\wp(u_i)}$. Again, since this condition is independent on the starting point, it follows that a generalized Poncelet porism holds as well.



\section{The ${\cal O}(4)$ spectral curve} \label{SEC:O4-multiplets}

\subsection{Representations}

In \cite{MM1} we began a detailed analysis of the ${\cal O}(4)$ multiplets and the elliptic curves associated to them. Before taking this analysis further, we review here its main points.

${\cal O}(4)$ multiplets can be written locally in either one of the following two generic forms
\begin{eqnarray}
\eta^{(4)}(\zeta) &\hspace{-7pt} = \hspace{-7pt}& \frac{\bar{z}}{\zeta^2}+\frac{\bar{v}}{\zeta}+x-v\zeta+z\zeta^2 \nonumber \\
&\hspace{-7pt} = \hspace{-7pt}&  \frac{\rho}{\zeta^2} \frac{(\zeta-\alpha)(\bar{\alpha}\zeta+1)}{1+|\alpha|^2}\frac{(\zeta-\beta)(\bar{\beta}\zeta+1)}{1+|\beta|^2} \label{eta4}
\end{eqnarray}
As discussed in \cite{MM2}, the antipodal conjugation-related reality constraint that they satisfy is preserved by a projective $SU(2)$ group of automorphisms of the Riemann sphere which act through birational transformations on the inhomogeneous coordinate $\zeta$. The projective component generates real scaling transformations while the $SU(2)$ component gives the multiplet a valuable rotational structure. Specifically, with respect to the induced action of the $SU(2)$ group, the polynomial coefficients form a spin-$2$ multiplet. Alternatively, this $SU(2)$ can be viewed as the group of isometries of the Riemann sphere endowed with the round metric of Fubini and Study. The roots $\alpha$, $\beta$ and their antipodal conjugates form on the sphere a constellation which rotates  rigidly under isometric transformations. The Fubini-Study distance between $\alpha$ and $\beta$ is an invariant, and is given explicitly by
\begin{equation}
\delta_{\alpha\beta} = 2 \arccos k_{\alpha\beta} = 2 \arcsin k'_{\alpha\beta} \label{FS}
\end{equation}
where the chordal distance and radius $k_{\alpha\beta}$ and $k'_{\alpha\beta}$ are expressed in terms of the roots as follows
\begin{equation}
k_{\alpha\beta} = \frac{|1+\bar{\alpha}\beta|}{\sqrt{(1+|\alpha|^2)(1+|\beta|^2)}} \hspace{30pt} {\rm and} \hspace{30pt}  k'_{\alpha\beta} = \frac{|\alpha-\beta|}{\sqrt{(1+|\alpha|^2)(1+|\beta|^2)}} \label{chordal_dist}
\end{equation}
Note that $k^2_{\alpha\beta} + k_{\alpha\beta}^{\prime2}  = 1$ and thus $0 < k_{\alpha\beta}, k_{\alpha\beta}' < 1$. When there is no risk of confusion, we denote them simply as $k$ and $k'$.

To each ${\cal O}(4)$ multiplet we associate the following quartic curve
\begin{equation}
\eta^2 = \zeta^2\eta^{(4)}(\zeta) \label{Majorana_n.f.}
\end{equation}
where $\eta$ is the second inhomogeneous coordinate on the complex plane. This is an algebraic curve of genus $1$ and a double cover of the real projective plane, $\mathbb{RP}^2$. We refer to this representation of the curve as the Majorana normal form. For a justification of this terminology, see \cite{MM2}. The ${\cal O}(4)$ curve can be birationally mapped to either the Legendre or the Weierstrass normal forms. In practice, this can be accomplished for instance with the two successive birational transformations given by\footnote{We use for cross-ratios the definition $\displaystyle{[z_1,z_2,z_3,z_4] = \frac{(z_1-z_3)(z_2-z_4)}{(z_1-z_4)(z_2-z_3)}}$.}
\begin{equation}
[\zeta,-\frac{1}{\bar{\alpha}},\ \alpha,\ \beta] = \nu = \frac{X-e_3}{e_1-e_3} \label{zeta->nu->X}
\end{equation}
where the Weierstrass roots are defined in terms of $k=k_{\alpha\beta}$ as follows
\begin{equation}
e_1 = -\frac{\rho}{3} (k^2-2) \qquad\quad e_2 = \frac{\rho}{3} (2k^2-1) \qquad\quad e_3 = -\frac{\rho}{3}(k^2+1) \label{Weierstrass_roots}
\end{equation}
The abelian differential form on the curve transforms accordingly
\begin{equation}
\varpi = \frac{d\zeta}{2\zeta\sqrt{\eta^{(4)}}} = \frac{d\nu}{2\sqrt{\rho\,\nu(1-\nu)(\nu-k^2)}} = \frac{dX}{2\sqrt{X^3 - g_2X - g_3}} \label{abelian-diff}
\end{equation}
Remarkably, the Weierstrass coefficients turn out to have explicit expressions in terms of the Majorana coefficients, {\it i.e.},
\begin{eqnarray}
g_2 &\hspace{-7pt} = \hspace{-7pt}& 4 |z|^2 + |v|^2 + \frac{1}{3} x^2  \label{g2_Majorana} \\
g_3 &\hspace{-7pt} = \hspace{-7pt}&  \frac{8}{3} |z|^2 x - \frac{1}{3} |v|^2 x - \frac{2}{27} x^3 - z \bar{v}^2 - \bar{z} v^2 \label{g3_Majorana}
\end{eqnarray}
Note that all the Legendre and Weierstrass moduli - the real scale $\rho$, the elliptic modulus $k$ and the Weierstrass coefficients $g_2$ and $g_3$ - are $SU(2)$-invariant.

The elliptic lattice $\Lambda = \mathbb{Z}\cdot 2\omega + \mathbb{Z}\cdot 2\omega'$ is generated by the loop integrals of the abelian differential form $\varpi$ over the canonical cycles of the ${\cal O}(4)$ curve. We have
\begin{equation}
\omega = \frac{K(k)}{\sqrt{\rho}} \qquad \mbox{and} \qquad \omega' = \frac{iK(k')}{\sqrt{\rho}}  \label{half-periods}
\end{equation}
As $0<k,k'<1$, both complete elliptic integrals $K(k)$ and $K(k')$ are real, and so the lattice is orthogonal.

\subsection{The ${\cal O}(4)$ Weierstrass cubic is a Cayley cubic} \label{SEC:Cayley-pencil}

If we define in place of the Majorana coefficients the related real variables
\begin{equation}
x_{\pm} = \frac{x\pm6|z|}{3} \qquad\qquad v_+ = \mbox{Im}\frac{v}{\sqrt{z}} \qquad\qquad v_- = \mbox{Re}\frac{v}{\sqrt{z}} \label{xpm_vpm}
\end{equation}
then in terms of these, the expressions (\ref{g2_Majorana}) and (\ref{g3_Majorana}) can be rewritten as follows
\begin{eqnarray}
g_2 &\hspace{-7pt} = \hspace{-7pt}& \ x_+^2+x_+x_-+x_-^2 + \frac{1}{4}(x_+\!-x_-)(v_-^2+v_+^2) \label{g2}   \\ 
g_3 &\hspace{-7pt} = \hspace{-7pt}& -(x_++x_-)x_+x_- - \frac{1}{4}(x_+\!-x_-)(x_+v_-^2+x_-v_+^2) \label{g3}
\end{eqnarray}
This form of the Weierstrass coefficients facilitates two key observations. First, we note that the four points with $(X,Y)$-coordinates 
\begin{equation}
\begin{array}{lcr}
\displaystyle{(x_-,v_-(x_+\!-x_-)/2)} &\quad& \displaystyle{(x_+,iv_+(x_+\!-x_-)/2)} \\ [9pt]
\displaystyle{(x_-,v_-(x_-\!-x_+)/2)} &\quad& \displaystyle{(x_+,iv_+(x_-\!-x_+)/2)}
\end{array} \label{four_points}
\end{equation}
are points on the ${\cal O}(4)$ curve in the Weierstrass representation, {\it i.e.} they satisfy the equation
\begin{equation}
Y^2 = X^3 - g_2X - g_3 \label{Weierstrass}
\end{equation}
This can be checked by direct substitution. The pairs of points on each column in (\ref{four_points}) are conjugated under the elliptic involution. One can also show that \cite{MM1}
\begin{equation}
e_3 < x_- < e_2 < x_+ < e_1 < e_0 = \infty \label{ineq_W}
\end{equation}
Secondly, we note that we can write the Weierstrass cubic as a determinant, {\it i.e.},
\begin{equation}
X^3 - g_2X - g_3 = 
\begin{array}{|ccc|}
X-x_+            & \sqrt{|z|} v_+ & 0 \\ [5pt]
\sqrt{|z|} v_+ & X+x_++x_-   & \sqrt{|z|} v_- \\ [5pt]
0 & \sqrt{|z|} v_- & X-x_-
\end{array}
\end{equation}
We give this fact the following interpretation: \\

\noindent {\it The Weierstrass cubic curve (\ref{Weierstrass}) associated to the ${\cal O}(4)$ multiplet is a Cayley cubic, {\it i.e.}},
\begin{equation}
Y^2 = \det ({\cal A}+X{\cal B})
\end{equation}
{\it for the pencil generated by the two plane conics with defining real-valued matrices}
\begin{equation}
{\cal A} = 
\left(\!\!\!
\begin{array}{ccc}
-x_+            & \sqrt{|z|} v_+ & 0 \\ [5pt]
\sqrt{|z|} v_+ & x_++x_-   & \sqrt{|z|} v_- \\ [5pt]
0 & \sqrt{|z|} v_- & -x_-
\end{array}
\!\!\!\right)
\hspace{20pt} \mbox{and} \hspace{20pt}
{\cal B} = \mathbb{I}_{3\times3} 
\label{Cayley-matrices}
\end{equation}

\subsection{Quantum spin coherent states} \label{SEC:spin-coh-st}

The relationship between ${\cal O}(2j)$ multiplets and spin coherent states has been extensively discussed in \cite{MM2}. Here, we review briefly some aspects of this correspondence that are going to play a role in our discussion.

The quantum states of a particle with spin $j$ are described in the spin coherent representation  by wave functions which are sections of ${\cal O}(2j)$ bundles over the so-called Bloch sphere, that is, they are polynomials of degree $j$ in the inhomogeneous coordinate on the sphere. Intuitively, such a state appears as a set of $2j$ elementary "spins 1/2" with the origins at the center of the Bloch sphere, pointing out in the directions marked by a constellation of $2j$ dots on the surface of the sphere corresponding to the roots of the wave function polynomial. A spin state is coherent when all elementary spins point in the same direction and is real when all elementary spins come in oppositely oriented pairs. Mathematically, the structure of the latter type of states is identical to that of the ${\cal O}(2j)$ multiplets discussed here.  The root factorization of the wave function polynomials can be viewed as a decomposition into spin-$1/2$ coherent states.

Spin-$1/2$ coherent states are in one-to-one correspondence to points on the Bloch sphere. Given a  spin-$1/2$ quantum system and $\zeta \in \mathbb{C} \cup \{\infty\} \simeq S^2$, the corresponding spin-$1/2$ coherent state is defined in terms of the standard orthonormal spin-up + spin-down basis in the Hilbert space of states of the system by the following linear superposition
\begin{equation}
|\zeta\rangle = \frac{1}{\sqrt{1+|\zeta|^2}} (\,|\!\downarrow\, \rangle + \zeta|\!\uparrow\, \rangle\,)
\end{equation}
The overlap between two spin-1/2 coherent states corresponding to $\alpha, \beta \in \mathbb{C} \cup \{\infty\}$ is
\begin{equation}
\langle \alpha | \beta \rangle = \frac{1+\bar{\alpha}\beta}{\sqrt{(1+|\alpha|^2)(1+|\beta|^2)}} 
\end{equation}
In particular, this formula implies that the overlap between states corresponding to antipodally opposite points is zero. The properties of quantum spin-$1/2$ coherent states are especially suited for use in spherical geometry, a feature that we will fully exploit later on. Thus, their norms are related to the Fubini-Study geodesic distance on the sphere between $\alpha$ and $\beta$. Specifically, with the definitions (\ref{chordal_dist}), one has
\begin{equation}
|\langle \alpha | \beta \rangle| = k_{\alpha\beta} \qquad\mbox{and}\qquad |\langle -\frac{1}{\bar{\alpha}} | \beta \rangle| = k'_{\alpha\beta} \label{coh_norms}
\end{equation}
On the other hand, the phases of cyclic sequences of spin-$1/2$ coherent states have an area interpretation, namely,
\begin{equation}
\langle \alpha_1 | \alpha_2 \rangle \langle \alpha_2 | \alpha_3 \rangle \cdots \langle \alpha_{n-1} | \alpha_n \rangle \langle \alpha_n | \alpha_1 \rangle = k_{\alpha_1\alpha_2} k_{\alpha_2\alpha_3} \cdots k_{\alpha_{n-1}\alpha_n} k_{\alpha_n\alpha_1} e^{i A_{\rm polygon}/2} \label{cyclic_coh}
\end{equation}
where $A_{\rm polygon}$ is the area of the spherical polygon with vertices at the points $\alpha_1$, \dots, $\alpha_n$. The factor $1/2$ in front of the area renders the ambiguity in the choice of what one means by the "inside" and the "outside" of the polygon irrelevant. 
\begin{figure}[htb]
\centering
\scalebox{0.5}{\includegraphics{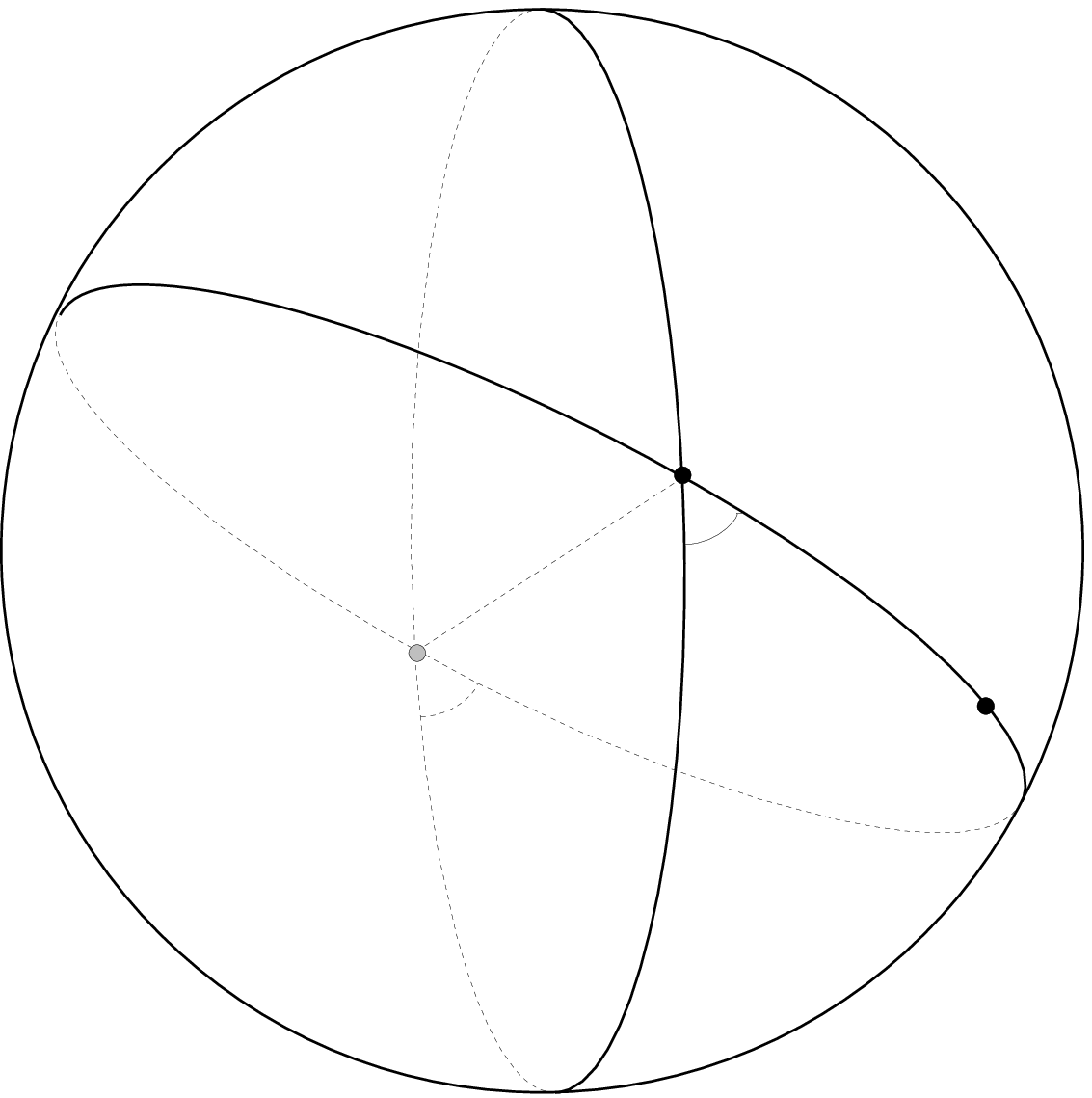}}
\put(-85,167){N}
\put(-84,-11){S}
\put(-58,74){$\phi_{\rm lune}$}
\put(-56,96){$\alpha$}
\put(-121,54){$\displaystyle{-\frac{1}{\bar{\alpha}}}$}
\put(-13,64){$\beta$}
\caption{}
\label{loone}
\end{figure} 
For later reference, let us also note that one can use equation (\ref{cyclic_coh}) to show that
\begin{equation}
\langle -\frac{1}{\bar{\alpha}} | \beta \rangle \langle \beta | \alpha \rangle = k'_{\alpha\beta} k_{\alpha\beta}\,  e^{i \phi_{\rm lune}} \label{la-luna}
\end{equation}
where $\phi_{\rm lune}$ is the dihedral angle of the {\it lune} cut on the sphere by the two geodesic circles that pass through $\beta$ and $\alpha$, respectively the South pole and $\alpha$, equal to half the area of the lune, see \figurename~\ref{loone}.

\subsection{The Jacobian picture} \label{SEC:Jacobian-pic}

One can check that, for any $\zeta \in \mathbb{C} \cup \{\infty\}$, one has
\begin{equation}
(\bar{X}_{\zeta}-e_2)(X_{-1/\bar{\zeta}}-e_2) = (e_1-e_2)(e_3-e_2) \label{dada}
\end{equation}
We use the notation $X_{\zeta}$ for the image of $\zeta$ through the birational map (\ref{zeta->nu->X}). Equation (\ref{dada}) implies that it is possible to choose the ambiguous signs of the $Y$-coordinates of the curve points with $X$-coordinates $\bar{X}_{\zeta}$ and $X_{-1/\bar{\zeta}}$ such that 
\begin{equation}
\begin{array}{|cll|}
\ \ \ 1 \ \ \   	& \bar{X}_{\zeta} 		& \bar{Y}_{\zeta}  \\ [5pt]
1 		& X_{-1/\bar{\zeta}} 	& Y_{-1/\bar{\zeta}}  \\ [5pt]
1 		& e_2 				& 0
\end{array}
= 0
\end{equation}
which is just the {\it colinearity condition} for the three points $(\bar{X}_{\zeta}, \bar{Y}_{\zeta})$, $(X_{-1/\bar{\zeta}}, Y_{-1/\bar{\zeta}})$ and $(e_2, 0)$. Moreover, this allows one to choose the corresponding points (through the Abel-Jacobi map) on the Jacobian variety such that
\begin{equation}
\bar{u}_{\zeta}+u_{-1/\bar{\zeta}} = \omega_2 \label{u+bu}
\end{equation}
This relation expresses the action of the antipodal conjugation-induced real structure on the Jacobian of the curve. 

By a straightforward calculation one can show that
\begin{eqnarray}
\cn (\sqrt{\rho}\, u_{\zeta}) &\hspace{-7pt} = \hspace{-7pt}& \sqrt{\frac{\alpha-\beta}{1+\bar{\alpha}\beta} \frac{1+\bar{\alpha}\zeta}{\alpha-\zeta}} \\
\dn (\sqrt{\rho}\, u_{\zeta}) &\hspace{-7pt} = \hspace{-7pt}& \sqrt{\frac{\alpha-\beta}{1+\bar{\beta}\beta} \frac{1+\bar{\beta}\zeta}{\alpha-\zeta}}
\end{eqnarray}
where $\cn$ and $\dn$ are the usual Jacobi elliptic functions. A similar expression holds for $\sn (\sqrt{\rho}\, u_{\zeta})$. These formulas are unsatisfactory for a number of reasons, chief among them being the fact that the roots $\alpha$ and $\beta$ do not appear on the same footing. We clearly need a different perspective. The crucial observation is contained in the following result
\begin{eqnarray}
\cn[\sqrt{\rho}\,(u_{\zeta} \pm \bar{u}_{\zeta})] 
&\hspace{-7pt} = \hspace{-7pt}&  \frac{k'_{\alpha\beta}}{k_{\alpha\beta}} \frac{k_{\alpha\zeta}k'_{\alpha\zeta} \mp k_{\beta\zeta}k'_{\beta\zeta}}{k_{\alpha\zeta}^{\prime2} - k_{\beta\zeta}^{\prime2}} 
= \frac{\displaystyle \tan \frac{\delta_{\alpha\beta}}{2}}{\displaystyle \tan \frac{\delta_{\alpha\zeta}\pm\delta_{\beta\zeta}}{2}} \label{cn(u+bu)} \\
\dn[\sqrt{\rho}\,(u_{\zeta} \pm \bar{u}_{\zeta})]  
&\hspace{-7pt} = \hspace{-7pt}&  k'_{\alpha\beta} \frac{k_{\alpha\zeta}k'_{\beta\zeta} \mp k_{\beta\zeta}k'_{\alpha\zeta}}{k_{\alpha\zeta}^2 - k_{\beta\zeta}^2} 
\, = \frac{\displaystyle \sin \frac{\delta_{\alpha\beta}}{2}}{\displaystyle \sin \frac{\delta_{\beta\zeta}\pm\delta_{\alpha\zeta}}{2}} \label{dn(u+bu)}
\end{eqnarray}
Incidentally, note that the these are the same type of trigonometric ratios that appear in the Napier and Delambre analogies of spherical trigonometry. The first equalites in (\ref{cn(u+bu)}) and (\ref{dn(u+bu)}) follow from applying the addition formulas (\ref{sncndn_add_thms}) for the Jacobi elliptic functions $\cn$ and $\dn$. We use that
\begin{equation}
\sn (\sqrt{\rho}\,u) = \sqrt{\frac{e_1-e_3}{X-e_3}} 
\qquad\quad
\cn (\sqrt{\rho}\,u) = \sqrt{\frac{X-e_1}{X-e_3}} 
\qquad\quad 
\dn (\sqrt{\rho}\,u) = \sqrt{\frac{X-e_2}{X-e_3}}
\end{equation}
with $X$ being the image of $u$ through the Weierstrass map, {\it i.e.}, $X=\wp(u)$. Despite the simple form of the outcome, the calculation is quite entangled and laborious if approached frontally. We managed to simplify and streamline it significantly by resorting to the spin coherent state techniques discussed above. First, observe that we have the following cross-ratio expressions
\begin{eqnarray}
\frac{X_{\zeta}-e_1}{e_3-e_1} &\hspace{-7pt} = \hspace{-7pt}& [\beta,-\frac{1}{\bar{\alpha}},\ \alpha,\ \zeta] \\ [2pt]
\frac{X_{\zeta}-e_2}{e_3-e_2} &\hspace{-7pt} = \hspace{-7pt}& [\beta,-\frac{1}{\bar{\beta}},\ \alpha,\ \zeta] \\
\frac{X_{\zeta}-e_3}{e_1-e_3} &\hspace{-7pt} = \hspace{-7pt}& [\zeta,-\frac{1}{\bar{\alpha}},\ \alpha,\ \beta] \label{X-e_3=cr}
\end{eqnarray}
Equation (\ref{X-e_3=cr}) is just a copy of (\ref{zeta->nu->X}); the preceding two follow from this one. The second observation is that cross-ratios can be expressed in terms of  spin-1/2 coherent states as follows
\begin{equation}
[z_1,z_2,z_3,z_4] = \frac{\displaystyle\langle-\frac{1}{\bar{z}_1}|z_3\rangle\langle-\frac{1}{\bar{z}_2}|z_4\rangle}{\displaystyle\langle-\frac{1}{\bar{z}_1}|z_4\rangle\langle-\frac{1}{\bar{z}_2}|z_3\rangle}
\end{equation}
Together, these relations allow one to cast the $\cn$ and $\dn$ addition formulas entirely in terms of spin-1/2 coherent states. The $k$ and $k'$ expressions emerge from the coherent state picture by means of the norm relations (\ref{coh_norms}). The second equalities in (\ref{cn(u+bu)}) and (\ref{dn(u+bu)}) follow by using further the relations (\ref{FS}) and some trigonometry.

For any $\zeta \in \mathbb{C} \cup \{\infty\}$, let us define 
\begin{equation}
u_{\zeta}^{\pm} = u_{\zeta} \pm u_{-1/\bar{\zeta}} \label{u^pm}
\end{equation}
{\it i.e.}, the "real" and "imaginary" parts of $u_{\zeta}$ with respect to the real structure induced by the antipodal conjugation on the sphere. Based on the equations (\ref{u+bu}), (\ref{cn(u+bu)}), (\ref{dn(u+bu)}) and the half-period addition formula $\sn [v \pm (K+iK')] = \pm\, \dn v \div  k\,\cn v$ we obtain
\begin{equation}
\sn(\sqrt{\rho}\, u_{\zeta}^{\pm}) =  \sec \frac{\delta_{\alpha\zeta}\mp\delta_{\beta\zeta}}{2} \label{sn(u^pm)}
\end{equation}
If we resort instead to the addition formula $\sn [v \pm K] = \pm\, \cn v \div \dn v$, we obtain
\begin{equation}
\sn[\sqrt{\rho}\,( u_{\zeta}^{\pm}\! - \omega')] = \frac{\displaystyle \cos \frac{\delta_{\alpha\zeta}\mp\delta_{\beta\zeta}}{2}}{\displaystyle \cos \frac{\delta_{\alpha\beta}}{2}} \label{sn(u^pm-om3)}
\end{equation}
We use here the conventional notations $K=K(k)$ and $K'=K(k')$ for the complete elliptic integrals of the first kind of complementary moduli. Remember now that we work on the 2-sphere with the antipodal points identified. This means essentially that we can always consider that the points $\alpha$, $\beta$ and $\zeta$ are on the {\it same hemisphere} of $S^2$. They determine a spherical triangle with vertices at $\alpha$, $\beta$, $\zeta$ and sides $\delta_{\alpha\zeta}$, $\delta_{\beta\zeta}$, $\delta_{\alpha\beta}$, which, for this reason, has the following properties: $1)$ $\delta_{\alpha\zeta}$, $\delta_{\beta\zeta}$, $\delta_{\alpha\beta}$ $\in [0,\pi]$, meaning the triangle is convex, which further implies that the usual triangle inequalities hold, {\it i.e.}, $\delta_{\alpha\zeta}+\delta_{\beta\zeta} \geq \delta_{\alpha\beta}$, {\it etc.} and $2)$ $\delta_{\alpha\zeta}+\delta_{\beta\zeta}+\delta_{\alpha\beta} \leq 2\pi$. Based on these inequalities being satisfied one determines that both equations (\ref{sn(u^pm)}) (that is, with both sets of signs considered) and the equation (\ref{sn(u^pm-om3)}) with the upper set of signs are $\geq 1$, whereas the equation (\ref{sn(u^pm-om3)}) with the lower set of signs is $\leq 1$ and $\geq -1$. It seems then natural to set this latter equation equal to the sine of an angle, let us call it $\sin D_{\zeta}$. In the light of (\ref{sine_thm}) we find it convenient to write this {\it definition} in the form
\begin{equation}
\frac{\displaystyle{\sin\frac{\pi-\delta_{\alpha\zeta}-\delta_{\beta\zeta}}{2}}}{\sin D_{\zeta}\rule{0pt}{12pt}} = k
\end{equation}
Inverting the lower equation (\ref{sn(u^pm-om3)}) on a fundamental domain yields
\begin{equation}
u_{\zeta}^{-} = \frac{1}{\sqrt{\rho}}F(\sin D_{\zeta},k) + \omega' \label{u_zeta^-}
\end{equation}
with $F(\cdot,k)$ an incomplete Legendre elliptic integral of the first kind. This provides us with a very explicit expression for $u_{\zeta}^-$, with a clearly resolved complex structure: the first term in the r.h.s. of (\ref{u_zeta^-}) is real, the second one is a purely imaginary constant shift.

We end this section with yet another important observation. We found that it is possible to choose the ambiguous signs of $Y_0$, $Y_{\infty}$ and $y_{\pm}$ corresponding on the Weierstrass curve (\ref{Weierstrass}) to $X_0$, $X_{\infty}$ and $x_{\pm}$, such that
\begin{equation}
\begin{array}{|ccc|}
 \  \ 1  \  \  	& X_{\infty} 			& \ Y_{\infty} \  \\ [5pt]
1 		& \hspace{-3pt} X_{0} 	& \hspace{-12pt}  \pm Y_{0}  \\ [5pt]
1 		& x_{\pm} 				& \hspace{-8pt} -y_{\pm}
\end{array}
= 0 \label{col-infty}
\end{equation}
This can be verified for instance by expressing everything in terms of the roots $\alpha$, $\beta$, their complex conjugates and the scale $\rho$, by means of the equations (\ref{zeta->nu->X}), (\ref{Weierstrass_roots}) and (\ref{chordal_dist}).
The equation (\ref{col-infty}) is a colinearity condition. By comparing the corresponding equation on the Jacobian to equation (\ref{u^pm}) with $\zeta = \infty$, we infer immediately that
\begin{equation}
\wp(u_{\infty}^{\pm}) = x_{\pm}
\end{equation}
{\it i.e.}, the four points $\pm u_{\infty}^+$ and $\pm u_{\infty}^-$ from the Jacobian are mapped by the inverse Abel-Jacobi map to the four points (\ref{four_points}) on the Weierstrass curve, with the $X$-coordinates equal to $x_+$ respectively $x_-$.

\subsection{${\cal O}(4)$ elliptic integrals} \label{integrals} \label{SEC:key-O4-int}

Generalized Legendre transform constructions based on ${\cal O}(4)$ multiplets oftentimes involve evaluating contour integrals of the type
\begin{equation}
{\cal I}_m = \int_{\Gamma} \frac{d\zeta}{\zeta} \frac{\zeta^m}{2\sqrt{\eta^{(4)}}} \label{I_m}
\end{equation}
with $\Gamma$ an integration contour which may be either open or closed, depending on context, and $m$ an integer taking values from $-2$ to $2$. In fact, it suffices to consider only $m = 0,1,2$, since the integrals corresponding to $m$ and $-m$ are complex conjugated to each other, modulo a shift. More precisely,
\begin{equation}
{\cal I}_{-m} = (-)^m \bar{\cal I}_m \pm 2\pi i \frac{m\bar{\beta}^{m-1}}{\sqrt{z}} \label{I_-m}
\end{equation}
This can be seen by changing in (\ref{I_m}) the integration variable $\zeta$ to $-1/\bar{\zeta}$ and then deforming the resulting contour back to the original one; in the process, one picks up a residue, which accounts for the shift term. Shifts are usually discarded by means of a doubling trick: we can always choose two contours, one which gives a $+$ and one which gives a $-$ in (\ref{I_-m}); by summing the two contributions up, the residue terms will mutually cancel. 

By "evaluating" these contour integrals we mean of course reducing them to standard elliptic integrals. For various reasons, we are particularly interested in obtaining as explicit a dependence on the Majorana coefficients of $\eta^{(4)}$ as possible. As it turns out, the Weierstrass framework is best suited to this purpose. Hence the first step of our approach is to transform the integrals from what we refer to as the Majorana picture to the Weierstrass picture by means of the birational transformation (\ref{zeta->nu->X}), which, with the help  of the notation that we have introduced at the begining of section \ref{SEC:Jacobian-pic}, can be conveniently written in the form
\begin{equation}
\zeta = \beta \frac{X-X_0}{\ X-X_{\infty}}
\end{equation}
The abelian differential that plays the role of integration measure transforms according to (\ref{abelian-diff}) as follows
\begin{equation}
\frac{d\zeta}{\zeta}\frac{1}{2\sqrt{\eta^{(4)}}} = \frac{dX}{2Y} \label{measure-MW}
\end{equation}
Once an integral is expressed completely in terms of Weierstrass variables, we follow the standard procedure in evaluating elliptic integrals, see {\it e.g.} \cite{MR698780}: we expand the rational coefficient of the measure (\ref{measure-MW}) into partial fractions centered on $X_{\infty}$ and then use formulas (\ref{Wei1}) through (\ref{Wei3}) to express each resulting term in terms of Weierstrass elliptic functions. That is of course not possible to do directly for the ${\cal I}_2$ integral, as the partial fraction expansion yields in that case a term proportional to $Y_{\infty}^2/(X-X_{\infty})^2$. One handles this by noticing that
\begin{equation}
\left(\!\frac{Y_{\infty}}{X-X_{\infty}}\!\right)^2 = \frac{1}{2}\!\left( X-X_{\infty} - \frac{3X_{\infty}^2-g_2}{X-X_{\infty}}\right) - Y\frac{d}{dX} \!\left(\!\frac{Y}{X-X_{\infty}}\!\right) \label{parfrac}
\end{equation}
The last term in (\ref{parfrac}) leads eventually to a total derivative which can be easily integrated. The other ones lead directly to elliptic integrals of the three kinds, just as in the other cases.

The outcome at this stage can be simplified by using that
\begin{eqnarray}
X_{\infty} &\hspace{-7pt} = \hspace{-7pt}& \frac{x}{3} - \beta v + 2\beta^2 z \\ [1pt]
\frac{3X_{\infty}^2-g_2}{Y_{\infty}} &\hspace{-7pt} = \hspace{-7pt}& \frac{v}{\sqrt{z}} - 4\beta\sqrt{z} \\
\frac{X_{\infty}-X_0}{Y_{\infty}} &\hspace{-7pt} = \hspace{-7pt}& -\frac{1}{\beta\sqrt{z}} 
\end{eqnarray}
and that
\begin{equation}
2\zeta(u_{\infty}) = \zeta(u_{\infty}^+) + \zeta(u_{\infty}^-) + 2\beta\sqrt{z}
\end{equation}
The first three identities can be verified by expressing everything in terms the Majorana roots and scale. The last one follows by applying succesively the doubling formula and then the addition theorem for the Weierstrass $\zeta$-function. Note also that (\ref{u^pm}) implies that $2u_{\infty} = u_{\infty}^+\!+u_{\infty}^-$.

In the end, we obtain
\begin{eqnarray}
{\cal I}_0 &\hspace{-7pt} = \hspace{-7pt}& u + {\cal C} \label{I0-inc} \\ [10pt]
{\cal I}_1 &\hspace{-7pt} = \hspace{-7pt}& - \frac{1}{2\sqrt{z}} \left[ \ln \frac{\sigma(u-u_{\infty})}{\sigma(u+u_{\infty})} + [\zeta(u_{\infty}^+)+\zeta(u_{\infty}^-)] u \right]  + {\cal C} \label{I1-inc} \\ [2pt]
{\cal I}_2 &\hspace{-7pt} = \hspace{-7pt}& -\frac{1}{4z} \left\{ \zeta(u-u_{\infty}) + \zeta(u+u_{\infty}) +  (x_+\!+x_-)u      
\phantom{ \frac{\sigma(u+u_{\infty})}{\sigma(u-u_{\infty})}} \right. \nonumber \\ [2pt]
&\hspace{-7pt}  \hspace{-7pt}& \left.  + \frac{v}{\sqrt{z}} \left[  \ln \frac{\sigma(u-u_{\infty})}{\sigma(u+u_{\infty})} + [\zeta(u_{\infty}^+)+\zeta(u_{\infty}^-)] u  \right] \right\} + {\cal C} \label{I2-inc}
\end{eqnarray}
where $u$ is related to $X$ as in equation (\ref{inv-AJ-map}). The corresponding complete integrals, obtained by integrating over the contours $\Gamma_i$  with $i = 1,2,3$ defined in the paragraph preceding equations (\ref{Wcei1}) through (\ref{Wcei3}), are
\begin{eqnarray}
{\cal I}_0^{(i)} &\hspace{-7pt} = \hspace{-7pt}& 2\omega_i \label{I0-c} \\ [8pt]
{\cal I}_1^{(i)} &\hspace{-7pt} = \hspace{-7pt}& \frac{1}{\sqrt{z}}[\pi_i(x_+)+\pi_i(x_-)]  \label{I1-c} \\
{\cal I}_2^{(i)} &\hspace{-7pt} = \hspace{-7pt}& -\frac{1}{2z} \left[ 2\eta_i + (x_+\!+x_-)\omega_i - \frac{v}{\sqrt{z}}[\pi_i(x_+)+\pi_i(x_-)]  \right] \label{I2-c}
\end{eqnarray}
To derive (\ref{I0-c}) - (\ref{I2-c}) from (\ref{I0-inc}) - (\ref{I2-inc}) we made use of the $\sigma$-function monodromy property as well as of a version of the $\zeta$-function addition theorem.

\section{ALE spaces of type $D_n$} \label{SEC:ALE-typeD}

\begin{figure}[hbt]
\centering
\scalebox{0.8}{\includegraphics{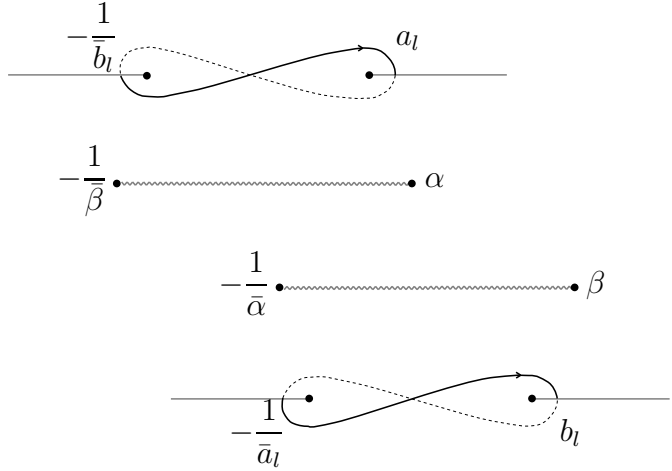}}
\put(-94,92){$\alpha$}
\put(-233,92){$\displaystyle{-\frac{1}{\bar{\beta}}}$}
\put(-33,52){$\beta$}
\put(-172,52){$\displaystyle{-\frac{1}{\bar{\alpha}}}$}
\put(-43,-4){$b_l$}
\put(-168,-4){$\displaystyle{-\frac{1}{\bar{a}_l}}$}
\put(-105,146){$a_l$}
\put(-230,146){$\displaystyle{-\frac{1}{\bar{b}_l}}$}
\caption{The two parts of the contour $\Gamma_l$.} \label{contour_gamma_l}
\end{figure}

The $F$-function that generates the asymptotically locally Euclidean (ALE) $D_n$ metric through the generalized Legendre transform construction is given, according to \cite{Chalmers:1998pu,Cherkis:1998xca,Cherkis:2003wk}, by
\begin{equation}
F = \oint_{\Gamma} \frac{d\zeta}{\zeta} \sqrt{\eta^{(4)}} - \sum_{l=1}^n \sum_{+,-} \frac{1}{2\pi i}\oint_{\Gamma_l} \frac{d\zeta}{\zeta} [\sqrt{\eta^{(4)}} \pm \chi^{(2)}_l] \ln [\sqrt{\eta^{(4)}} \pm \chi^{(2)}_l] \label{F_Dn}
\end{equation}
The parameters of the ${\cal O}(2)$-multiplets $\chi^{(2)}_l$, which transform as the components of a vector at rotations, do not coordinatize the ALE space but rather specify the positions of the monopoles. The contour $\Gamma$ winds around the canonical 2-cycles of  $\sqrt{\eta^{(4)}}$. The $n$ contours $\Gamma_l$ surround the roots $a_l$, $-1/\bar{a}_l$, $b_l$, $-1/\bar{b}_l$ of the deformed ${\cal O}(4)$ multiplets $\eta^{(4)} - (\chi^{(2)}_l)^2$ in the way depicted schematically in \figurename~\ref{contour_gamma_l}. 

The roots $a_l$, $-1/\bar{a}_l$, $b_l$, $-1/\bar{b}_l$ are obtained by solving for $\zeta$ the equation
\begin{equation}
\eta^{(4)}(\zeta) = \eta_l^{(2)}(\zeta){}^2
\end{equation}
This is an equation on $\mathbb{RP}^2$, the 2-sphere with antipodal points identified. We can get some insight into it by using the spin-$1/2$ coherent wave-function representation of section \ref{SEC:spin-coh-st}. In terms of it, the equation can be rewritten as follows
\begin{equation}
\rho\, \langle -\frac{1}{\bar{\zeta}} | \alpha \rangle \langle \alpha | \zeta \rangle \langle -\frac{1}{\bar{\zeta}} | \beta \rangle \langle \beta | \zeta \rangle = \left( \sigma_l   \langle -\frac{1}{\bar{\zeta}} | \gamma_l \rangle \langle \gamma_l | \zeta \rangle \right)^{\!2} 
\end{equation}
Then, based on the equations (\ref{coh_norms}) and (\ref{la-luna}), by taking the norm and, separately, comparing the phase factors on the two sides, one obtains
\begin{eqnarray}
 \rho\, \sin \delta_{\alpha\zeta} \sin \delta_{\beta\zeta} &\hspace{-7pt} = \hspace{-7pt}& \sigma_l^2 \sin^2 \!\delta_{\gamma_l\zeta} \label{nrm} \\ [5pt]
 \phi_{\widehat{\gamma_l\zeta\alpha}} + \phi_{\widehat{\gamma_l\zeta\beta}} &\hspace{-7pt} = \hspace{-7pt}& 2\pi k \quad\quad (k \in \mathbb{Z}) \label{dihedral}
\end{eqnarray}
where $ \phi_{\widehat{\gamma_l\zeta\alpha}}$ is the (oriented) angle formed by  the two geodesic circles that pass through $\gamma_l$ respectively $\alpha$ and intersect at $\zeta$; $\phi_{\widehat{\gamma_l\zeta\beta}}$ is defined similarly. Equation (\ref{dihedral}) means geometrically that $\gamma_l$ sits on the geodesic circle that bisects the angle formed by the two geodesic circles that pass through $\alpha$ respectively $\beta$ and intersect at $\zeta$. 
Unfortunately we do not yet possess a satisfactory understanding of the geometric picture behind these equations. But notice that if we think of them not as equations for $\zeta$ but for $\gamma_l$, or, in other words, if we formulate the problem in this way: given $\zeta$ fixed ($\alpha$ and $\beta$ are assumed fixed in either case), find $\gamma_i$ that leads to it, then a simple geometric picture emerges. In this case, equation (\ref{nrm}) can be easily solved to yield $\delta_{\gamma_l\zeta}$. Clearly, since we assume that $\delta_{\alpha\zeta}$, $\delta_{\beta\zeta}$, $\delta_{\gamma_l\zeta}$ $\in [0,\pi]$, one can have either no solution or two solutions (two supplementary angles), counting multiplicities. Notice that if $\sigma_l^2 \geq \rho$ then one always has two solutions. So let us assume there are two solutions. Arrange the sphere such that $\zeta$ and $-1/\bar{\zeta}$ lie on the North-South axis. Then the locus of $\gamma_l$ corresponding to a given pair of solutions for $\delta_{\gamma_l\zeta}$ is given by two circles parallel to the equator. But, as we stated above, the locus of solutions of equation (\ref{dihedral}) is the geodesic circle that bisects the spherical angle $\widehat{\alpha\zeta\beta}$ - a meridian, in our picture. The solutions for $\gamma_l$ lie at the intersection of the pair of paralel circles with this meridian. Note that these solutions come in antipodally-conjugated pairs, as objects that descend on $\mathbb{RP}^2$ should.

Denoting with $F_{\chi}$ the sum of $\chi$-deformed terms in (\ref{F_Dn}), then by commuting the derivatives with the integrals one obtains
\begin{eqnarray}
\frac{\partial F_{\chi}}{\partial x} 
&\hspace{-7pt} = \hspace{-7pt}& -\sum_{l=1}^n \frac{1}{2\pi i} \oint_{\Gamma_l} \frac{d\zeta}{\zeta} \frac{1}{2\sqrt{\eta^{(4)}}} \ln [\eta^{(4)} - (\chi^{(2)}_l)^2] \nonumber \\ [2pt]
&\hspace{-7pt} = \hspace{-7pt}& -\sum_{l=1}^n \int^{\ a_l}_{-1/\bar{b}_l} \hspace{-5pt}+ \int^{\ b_l}_{-1/\bar{a}_l} \frac{d\zeta}{\zeta} \frac{1}{2\sqrt{\eta^{(4)}}}
\end{eqnarray}
and
\begin{eqnarray}
\frac{\partial F_{\chi}}{\partial v} 
&\hspace{-7pt} = \hspace{-7pt}& \sum_{l=1}^n \frac{1}{2\pi i} \oint_{\Gamma_l} \frac{d\zeta}{\zeta} \frac{\zeta}{2\sqrt{\eta^{(4)}}} \ln [\eta^{(4)} - (\chi^{(2)}_l)^2] \nonumber \\ [2pt]
&\hspace{-7pt} =  \hspace{-7pt}& \sum_{l=1}^n \int^{\ a_l}_{-1/\bar{b}_l} \hspace{-5pt}+ \int^{\ b_l}_{-1/\bar{a}_l} \frac{d\zeta}{\zeta} \frac{\zeta}{2\sqrt{\eta^{(4)}}}
\end{eqnarray}
The logarithm can be dropped out of the integral at the expense of turning closed contours into open contours. We thus arrive at {\it incomplete} elliptic integrals of the type (\ref{I_m}), with $m=0,1$. The first integral in (\ref{F_Dn}) appears also in the Atiyah-Hitchin case \cite{MM1} and leads to {\it complete} elliptic integrals of the same type. Using the fundamental results of section \ref{integrals} we derive in a straightforward manner the following formulas
\begin{equation}
\frac{\partial F}{\partial x} = 2m\omega + 2m'\omega' - \sum_{l=1}^n (u_{a_l}^- +u_{b_l}^-) \label{F_x}
\end{equation}
and
\begin{eqnarray}
\frac{\partial F}{\partial v} &\hspace{-7pt} = \hspace{-7pt}& \frac{1}{2\sqrt{z}} \ln \frac{\sigma(2m\omega \!-\! u_{\infty})\sigma(2m'\omega' \!-\! u_{\infty})}{\sigma(2m\omega \!+\! u_{\infty})\sigma(2m'\omega' \!+\! u_{\infty})}  \prod_{l=1}^n \prod_{\zeta = a_l, b_l} \!\! \frac{\sigma(u_{\zeta} \!+\! u_{\infty})\sigma(u_{-1/\bar{\zeta}} \!-\! u_{\infty})}{\sigma(u_{\zeta} \!-\! u_{\infty})\sigma(u_{-1/\bar{\zeta}} \!+\! u_{\infty})} \nonumber \\
&\hspace{-7pt} + \hspace{-7pt}& \frac{1}{2\sqrt{z}} [\zeta(u_{\infty}^+)+\zeta(u_{\infty}^-)] \frac{\partial F}{\partial x} \label{F_v}
\end{eqnarray}
with $m, m' \in \mathbb{Z}$. Observe that equation (\ref{F_x}) determines the winding number $m'$ if we require that $F$ be real. In this case, since $x$ is real, the whole equation has to be real. From (\ref{u_zeta^-}) it is clear that the imaginary parts of both $u_{a_l}^-$ and $u_{b_l}^-$ are equal to $\omega'$. To cancel them, one needs to take $m' = n$.

Since the $n$ multiplets $\chi^{(2)}_l$ are spectators, the Legendre relations read 
\begin{eqnarray}
\frac{\partial F}{\partial v} &\hspace{-7pt} = \hspace{-7pt}& u \label{L-rel1} \\ [2pt]
\frac{\partial F}{\partial x} &\hspace{-7pt} = \hspace{-7pt}& 0 \label{L-rel2}
\end{eqnarray}
Together with equation (\ref{F_v}) they give
\begin{equation}
e^{2u\sqrt{z}} = \frac{\sigma(2m\omega \!-\! u_{\infty})\sigma(2m'\omega' \!-\! u_{\infty})}{\sigma(2m\omega \!+\! u_{\infty})\sigma(2m'\omega' \!+\! u_{\infty})}  \prod_{l=1}^n \prod_{\zeta = a_l, b_l} \!\!  \frac{\sigma(u_{\zeta} \!+\! u_{\infty})\sigma(u_{-1/\bar{\zeta}} \!-\! u_{\infty})}{\sigma(u_{\zeta} \!-\! u_{\infty})\sigma(u_{-1/\bar{\zeta}} \!+\! u_{\infty})}
\end{equation}
The expression on the r.h.s. is a meromorphic elliptic function in $u_{\infty}$, with zeros at $2m\omega$, $2m'\omega'$, $-u_{a_l}$, $u_{-1/\bar{a}_l}$, $-u_{b_l}$, $u_{-1/\bar{b}_l}$ for all values of $l$, and poles at the mirror points, of opposite sign.

On the other hand, equation (\ref{L-rel2}) together with the expression (\ref{F_x}) imply
\begin{equation}
\sum_{l=1}^n [F(\sin D_{a_l},k) + F(\sin D_{b_l},k)] = \mathbb{Z} \cdot 2K(k) \label{sumF=Z2K}
\end{equation}
where the angles $D_{a_l}$ and $D_{b_l}$ are defined by
\begin{equation}
\frac{\displaystyle{\sin\frac{\pi-\delta_{\alpha a_l}-\delta_{\beta a_l}}{2}}}{\sin D_{a_l}\rule{0pt}{12pt}} = 
\frac{\displaystyle{\sin\frac{\pi-\delta_{\alpha b_l}-\delta_{\beta b_l}}{2}}}{\sin D_{b_l}\rule{0pt}{12pt}} = k
\qquad \mbox{ for $l=1,\cdots,n$}
\end{equation}
We write these relations in this form in order to make the resemblance to Legendre's addition theorem expressed by equations (\ref{Legendre_add_thm}) and (\ref{sine_thm}) transparent. This is the Ercolani-Sinha-type constraint to which we alluded in the introduction.

In section \ref{SEC:Cayley-pencil} we have established that the ${\cal O}(4)$ curve has,  in the Weierstrass representation, a natural Cayley pencil structure for the two plane conics with defining matrices ${\cal A}$ and ${\cal B}$ given explicitly in terms of the parameters of the ${\cal O}(4)$ polynomial in (\ref{Cayley-matrices}). On the other hand, the equation (\ref{sumF=Z2K}) has precisely the form of an addition theorem of the type (\ref{gen-cycl-n}) on the Jacobian of the curve. In the light of section \ref{SEC:Poncelet-porism}, these facts suggest the following geometric quantization interpretation: \\

\noindent {\it The generalized Legendre relation (\ref{L-rel2}) takes in the case of the ALE spaces of type $D_n$ the form of a closure condition for a Poncelet polygon with vertices lying on the conic ${\cal B}$ and sides tangent to various conics of the pencil generated by ${\cal A}$ and ${\cal B}$ determined by the positions of the monopoles. Poncelet's porism ensures that this condition is not dependent on the starting point in the construction of the polygon.} \\





\section[APPENDIX: Weierstrass elliptic integrals]{APPENDIX: Weierstrass elliptic integrals\footnote{We reproduce this Appendix {\it verbatim} from \cite{MM1}, to which we also refer for further details.}}

In the Weierstrass theory the role of the incomplete elliptic integrals is played by 
\begin{eqnarray}
&& \phantom{+} \int \frac{dX}{2Y} = u + {\cal C} \label{Wei1} \\ [3pt]
&& -\int X\frac{dX}{2Y} = \zeta(u) + {\cal C} \label{Wei2} \\ [2pt]
&& -\int \frac{Y_0}{X-X_0} \frac{dX}{2Y} = \frac{1}{2} \ln \frac{\sigma(u+u_0)}{\sigma(u-u_0)} - u\,\zeta(u_0) + {\cal C} \label{pi_W} \label{Wei3}
\end{eqnarray}
where ${\cal C}$ is an indefinite integration constant, $(X,Y)$ and  $(X_0,Y_0)$ are points on the Weierstrass cubic $Y^2 = X^3-g_2X-g_3$, $u$ and $u_0$ are the corresponding points on the Jacobian variety, and $\sigma(u)$, $\zeta(u)$ are the Weierstrass sigma respectively zeta pseudo-elliptic functions. The expressions on the r.h.s. are obtained by substituting $X$ and $Y$ with the corresponding Weierstrass elliptic functions, {\it i.e.},
\begin{equation}
X = \wp(u;4g_2,4g_3) \hspace{30pt}  2Y = \wp'(u;4g_2,4g_3) \label{inv-AJ-map}
\end{equation}
The derivation of the first two expressions is fairly straightforward and standard. The derivation of the third one requires the use of a variant of the addition theorem of the Weierstrass zeta-function. 

The corresponding complete integrals are obtained by integrating in the complex $X$-plane along the closed countours $\Gamma_1$, surrounding the roots $e_2$ and $e_3$, $\Gamma_2$, surrounding the roots $e_3$ and $e_2$ and $\Gamma_3$, surrounding the roots $e_2$ and $e_1$, or, more precisely, on the Jacobian, from $u = \omega_2$ to $-\omega_3$, from $u = \omega_3$ to $-\omega_2$ and from $u = \omega_2$ to $-\omega_1$, respectively. We get
\begin{eqnarray}
&& \phantom{+} \oint_{\Gamma_i} \frac{dX}{2Y} = 2\omega_i \label{Wcei1} \\ [2pt]
&& -\oint_{\Gamma_i} X\frac{dX}{2Y} = 2\eta_i \label{Wcei2} \\ [2pt]
&& -\oint_{\Gamma_i} \frac{Y_0}{X-X_0} \frac{dX}{2Y}  = 2\
\begin{array}{|cc|}
u_0		& \omega_i \\
\zeta(u_0)	\!&\! \zeta(\omega_i) 
\end{array} 
\stackrel{\rm def}{=} 2\,\pi_i(X_0) \label{Wcei3}
\end{eqnarray}
where $u_0$ is the image of $(X_0,Y_0)$ through the Abel-Jacobi map and $i = 1,2,3$. Equation (\ref{Wcei3}) follows by way of the monodromy property of the Weierstrass sigma-function in the r.h.s. of (\ref{pi_W}). The notation $\pi_i(X_0)$ is not quite rigorous, a more appropriate one would be for instance $\pi_i(X_0,Y_0)$ or $\pi_i(u_0)$. We use it  nevertheless, but with the implicit {\it caveat} that it conceals a sign ambiguity. Clearly, only two out of three integrals of each set of integrals are independent, as $\omega_1+\omega_2+\omega_3 = 0$, $\eta_1+\eta_2+\eta_3 = 0$ and $\pi_1(X)+\pi_2(X)+\pi_3(X) = 0$.

\bibliographystyle{utphys}
\bibliography{article3}

\end{document}